\documentclass[11pt,a4paper]{article}
\usepackage[latin1]{inputenc}
\usepackage{amsmath}
\usepackage{amsfonts}
\usepackage{amssymb}
\usepackage{graphicx}
\usepackage{float}
\usepackage{stmaryrd}
\restylefloat{figure}

\newtheorem{satz}{Satz}[section]
\newtheorem{theorem}[satz]{Theorem}
\newtheorem{cor}[satz]{Corollary}
\newtheorem{defin}[satz]{Definition}
\newtheorem{lemma}[satz]{Lemma}

\newcommand{\abs}[1]{\left|{#1}\right|}
\newcommand{\norm}[1]{\left|\left|{#1}\right|\right|}
\newcommand{\rund}[1]{\left(#1\right)}
\newcommand{\spitz}[1]{\left\langle{#1}\right\rangle}
\newcommand{\schweif}[1]{\left\{#1\right\}}
\newcommand{\eckig}[1]{\left\lbrack{#1}\right\rbrack}

\def\zz{\mathbb{Z}}
\def\cz{\mathbb{C}}
\def\nz{{\rm I\kern-.20em N}}
\def\rz{{\rm I\kern-.20em R}}
\def\tz{\mathbb{T}}

\def\B{\begin{cal} B \end{cal}}
\def\C{\begin{cal} C \end{cal}}
\def\D{\begin{cal} D \end{cal}}
\def\H{\begin{cal} H \end{cal}}
\def\O{\begin{cal} O \end{cal}}
\def\S{\begin{cal} S \end{cal}}

\def\g{\frak g}
\def\a{\frak a}
\def\n{\frak n}
\def\k{\frak k}
\def\sl{\frak sl}
\def\m{\frak m}
\def\z{\frak z}

\def\vol{{\rm vol\ }}
\def\dist{{\rm dist}}
\def\Im{{\rm Im\ }}
\def\tr{{\rm tr}}
\def\Leb{{\rm Leb}}
\def\Aut{{\rm Aut}}
\def\cosh{{\rm cosh\ }}
\def\sinh{{\rm sinh\ }}
\def\Re{{\rm Re\ }}
\def\Ad{{\rm Ad}}

\def\bz{{\bf z} }
\def\bv{{\bf v} }
\def\be{{\bf e} }
\def\bb{{\bf b} }
\def\bc{{\bf c} }
\def\bu{{\bf u} }
\def\bw{{\bf w} }
\def\bs{{\bf s} }

\def\b0{{\bf 0} }
\def\bZ{{\bf Z} }
\def\bX{{\bf X} }

\def\eps{\varepsilon }
\def\bzeta{{\bf \zeta} }
\def\bEta{{\bf \eta} }
\def\bthet{{\bf \vartheta} }

\begin{document}

\vskip 1.0 true cm

\begin{center}
{\Huge \bf A Spanning Set for the space of Super Cusp forms}
\end{center}

\vskip 1.0 true cm

\begin{center}
Roland Knevel ,

Unité de Recherche en Mathématiques Luxembourg
\end{center}

\vskip 2.0 true cm

\section*{Mathematical Subject Classification}

11F55 (Primary) , 32C11 (Secondary) .

\section*{Keywords}

Automorphic and cusp forms, super symmetry, semisimple {\sc Lie} groups, partially hyperbolic flows, unbounded realization of a complex bounded symmetric domain.

\section*{Abstract}

Aim of this article is the construction of a spanning set for the space $sS_k(\Gamma)$ of super cusp forms on a complex bounded symmetric super domain $\B$ of rank $1$ with respect to a lattice $\Gamma$ . The main ingredients are a generalization of 
the {\sc Anosov} closing lemma for partially hyperbolic diffeomorphisms and an unbounded realization $\H$ of $\B$ , in particular {\sc Fourier} decomposition at the cusps of the quotient $\Gamma \backslash B$ mapped to $\infty$ via a partial {\sc Cayley}
transformation. The elements of the spanning set are in finite-to-one correspondence with closed geodesics of the body $\Gamma \backslash B$ of $\Gamma \backslash \B$ , the number of elements corresponding to a geodesic growing linearly with its 
length.

\section*{Introduction}

Automorphic and cusp forms on a complex bounded symmetric domain $B$ are already a well established field of research in mathematics. They play a fundamental role in representation theory of semisimple {\sc Lie} groups of Hermitian type, and they 
have applications to number theory, especially in the simplest case where $B$ is the unit disc in $\cz$ , biholomorphic to the upper half plane $H$ via a {\sc Cayley} transform, $G = SL(2, \rz)$ acting on $H$ via {\sc Möbius} transformations and
$\Gamma \sqsubset SL(2, \zz)$ of finite index. Aim of the present paper is to generalize an approach used by Tatyana {\sc Foth} and Svetlana {\sc Katok} in \cite{Foth} and \cite{Katok} for the construction of spanning sets for the space of cusp forms on a 
complex bounded symmetric domain $B$ of rank $1$~, which then by classification is (biholomorphic to) the unit ball of some $\cz^n$ , $n \in \nz$~, and a lattice $\Gamma \sqsubset G = \Aut_1(B)$ for sufficiently high weight $k$~. This is done in theorem
\ref{main} , which is the main theorem of this article, again for sufficiently large weight $k$ . \\

The new idea in \cite{Foth} and \cite{Katok} is to use the concept of a hyperbolic (or {\sc Anosov}) diffeomorphism resp. flow on a Riemannian manifold and an appropriate version of the {\sc Anosov} closing lemma. This concept originally comes from the 
theory of dynamical systems, see for example in \cite{KatHas} . Roughly speaking a flow $\rund{\varphi_t}_{t \in \rz}$ on a Riemannian manifold $M$ is called hyperbolic if there exists an orthogonal and $\rund{\varphi_t}_{t \in \rz}$-stable splitting 
$T M = T^+ \oplus T^- \oplus T^0$ of the tangent bundle $T M$ such that the differential of the flow $\rund{\varphi_t}_{t \in \rz}$ is uniformly expanding on $T^+$ , uniformly contracting on $T^-$ and isometric on $T^0$ , and finally $T^0$ is one-dimensional 
generated by $\partial_t \varphi_t$ . In this situation the {\sc Anosov} closing lemma says that given an 'almost' closed orbit of the flow $\rund{\varphi_t}_{t \in \rz}$ there exists a closed orbit 'nearby' . Indeed given a complex bounded symmetric domain 
$B$ of rank $1$~, $G = \Aut_1(B)$ is a semisimple {\sc Lie} group of real rank $1$ , and the root space decomposition of its {\sc Lie} algebra $\g$ with respect to a {\sc Cartan} subalgebra $\a \sqsubset \g$ shows that the geodesic flow 
$\rund{\varphi_t}_{t \in \rz}$ on the unit tangent bundle $S(B)$ , which is at the same time the left-invariant flow on $S(B)$ generated by $\a \simeq \rz$ , is hyperbolic. The final result in this direction is theorem \ref{Anosov} (i) . \\

For the super case first it is necessary to develop the theory of super automorphic resp. cusp forms, while the general theory of ($\zz_2$-) graded structures and super manifolds is already well established, see for example \cite{Const} . It has first 
been developed by F. A. {\sc Berezin} as a mathematical method for describing super symmetry in physics of elementary particles. However even for mathematicians the elegance within the theory of super manifolds is really amazing and satisfying.
Here I deal with a simple case of super manifolds, namely complex super domains. Roughly speaking a complex super domain $\B$ is an object which has $(n,r) \in \nz^2$ as super dimension and which has the characteristics:

\begin{itemize}
\item[(i)] it has a body $B = \B^\#$ being an ordinary domain in $\cz^n$ , 
\item[(ii)] the complex unital graded commutative algebra $\O(\B)$ of holomorphic super functions on $\B$ is (isomorphic to) $\O(B) \otimes \bigwedge\rund{\cz^r}$ , where $\bigwedge\rund{\cz^r}$ denotes the exterior algebra of $\cz^r$ . Furthermore 
$\O(\B)$ naturally embeds into the first two factors of the complex unital graded commutative algebra $\D(\B) \simeq \C^\infty(B)^\cz \otimes \bigwedge\rund{\cz^r} \boxtimes \bigwedge\rund{\cz^r} \simeq \C^\infty(B)^\cz \otimes \bigwedge\rund{\cz^{2 r}}$ of 
'smooth' super functions on $\B$ , where $\C^\infty(B)^\cz = \C^\infty(B, \cz)$ denotes the algebra of ordinary smooth functions with values in $\cz$ , which is at the same time the complexification of $\C^\infty(B)$ , and '$\boxtimes$' denotes the graded 
tensor product.
\end{itemize}

We see that for each pair $(B, r)$ where $B \subset \cz^n$ is an ordinary domain and $r \in \nz$ there exists exactly one $(n, r)$-dimensional complex super domain $\B$ of super dimension $(n, r)$ with body $B$ , and we denote it by $B^{| r}$ . Now let 
$\zeta_1, \dots, \zeta_n \in \cz^r$ denote the standard basis vectors of $\cz^r$ . Then they are the standard generators of $\bigwedge\rund{\cz^r}$ , and so we get the standard even (commuting) holomorphic coordinate functions 
$z_1, \dots, z_n \in \O(B) \hookrightarrow \O\rund{B^{ |r}}$ and odd (anticommuting) coordinate functions $\zeta_1, \dots, \zeta_r \in \bigwedge\rund{\cz^r} \hookrightarrow \O\rund{B^{ |r}}$~. So omitting the tensor products as there is no danger of confusion 
we can decompose every $f \in \O\rund{B^{ |r}}$ uniquely as

\[
f = \sum_{I \in \wp(r)} f_I \zeta^I   \, ,
\]

where $\wp(r)$ denotes the power set of $\{1, \dots, r\}$ , all $f_I \in \O(B)$ , $I \in \wp(r)$ , and $\zeta^I := \zeta_{i_1} \cdots \zeta_{i_s}$ for all $I = \schweif{i_1, \dots, i_s} \in \wp(r)$ , $i_1 < \dots < i_s$ . \\

$\D\rund{B^{ |r}}$ is a graded ${}^*$-algebra, and the graded involution

\[
\overline{\phantom{1}} : \D\rund{B^{ |r}} \rightarrow \D\rund{B^{ |r}}
\]

is uniquely defined by the rules

\begin{itemize}
\item[\{i\}] $\overline{\overline f} = f$ and $\overline{f h} = \overline h \, \, \overline f$ for all $f, h \in \D\rund{B^{ |r}}$ ,
\item[\{ii\}] $\overline{\phantom{1}}$ is $\cz$-antilinear, and restricted to $\C^\infty(B)$ it is just the identity,
\item[\{iii\}] $\overline{\zeta_i}$ is the $i$-th standard generator of $\bigwedge\rund{\cz^r} \hookrightarrow \D\rund{B^{ |r}}$ embedded as {\bf third} factor, where $\zeta_i$ denotes the $i$-th odd holomorphic standard coordinate on $B^{ |r}$ , which is the 
$i$-th standard generator of $\bigwedge\rund{\cz^r} \hookrightarrow \D\rund{B^{ |r}}$ embedded as {\bf second} factor, $i = 1, \dots, r$ .
\end{itemize}

With the help of this graded involution we are able to decompose every $f \in \D\rund{B^{ |r}}$ uniquely as

\[
f = \sum_{I, J \in \wp(r)} f_{I J} \zeta^I \, \overline \zeta^J \, ,
\]

where $f_{I J} \in \C^\infty(B)^\cz$ , $I, J \in \wp(r)$ , and $\overline \zeta^J := \overline{\zeta_{i_1}} \dots \overline{\zeta_{i_s}}$ for all $J = \schweif{j_1, \dots, j_s} \in \wp(r)$ , $j_1 < \dots < j_s$ . \\

For a discussion of super automorphic and super cusp forms we restrict ourselves to the case of the {\sc Lie} group $G := sS\rund{U(n, 1) \times U(r)}$ , $n \in \nz \setminus \{0\}$~, $r \in \nz$~, acting on the complex $(n,r)$-dimensional super unit ball 
$B^{| r}$ . So far there seems to be no classification of super complex bounded symmetric doimains although we know the basic examples, see for example in chapter IV of \cite{Borth} , which I follow here. The group $G$ is the 
body of the super {\sc Lie} group $SU(n, 1| r)$ studied in \cite{Borth} acting on $B^{| r}$ . The fact that an ordinary discrete subgroup (which means a sub super {\sc Lie} group of super dimension $(0, 0)$ ) of a super {\sc Lie} group is just an ordinary 
discrete subgroup of the body justifies our restriction to an ordinary {\sc Lie} group acting on $B^{| r}$ since purpose of this article is to study automorphic and cusp forms with respect to a lattice. In any case one can see the odd directions of the complex 
super domain $B^{ |r}$ already in $G$ since it is an almost direct product of the semisimple {\sc Lie} group $SU(n, 1)$ acting on the body $B$ and $U(r)$ acting on $\bigwedge\rund{\cz^r}$ . Observe that if $r > 0$ the full automorphism group of $B^{| r}$ , 
without any isometry condition, is never a super {\sc Lie} group since one can show that otherwise its super {\sc Lie} algebra would be the super {\sc Lie} algebra of integrable super vector fields on $B^{| r}$ , which has unfortunately infinite dimension. \\

Let me remark two striking facts:

\begin{itemize}
\item[(i)] the construction of our spanning set uses {\sc Fourier} decomposition exactly three times, which is not really surprising, since this corresponds to the three factors in the {\sc Iwasawa} decomposition $G = K A N$ .
\item[(ii)] super automorphic resp. cusp forms introduced this way are equivalent (but not one-to-one) to the notion of 'twisted' vector-valued automorphic resp. cusp forms.
\end{itemize}

{\it Acknowledgement:} Since the research presented in this article is partially based on my PhD thesis I would like to thank my doctoral advisor Harald {\sc Upmeier} for mentoring during my PhD but also Martin {\sc Schlichenmaier} and Martin {\sc Olbrich}
for their helpful comments.

\pagebreak

\section{The space of super cusp forms}

Let $n \in \nz \setminus \{0\}$ , $r \in \nz$ and

\begin{eqnarray*}
&& G := sS\rund{U(n, 1) \times U(r)} \\
&& \phantom{12} := \schweif{\left.\rund{\begin{array}{c|c}
g' & 0 \\ \hline
0 & E
\end{array}} \in U(n, 1) \times U(r) \, \right| \, \det g' = \det E}  \, ,
\end{eqnarray*}

which is a real $\rund{(n + 1)^2 + r^2 - 1}$-dimensional {\sc Lie} group. Let $\B := B^{| r}$ , where

\[
B := \schweif{\left. \bz \in \cz^n \, \right| \, \bz ^* \bz < 1} \subset \cz^n
\]

denotes the usual unit ball, with even coordinate functions $z_1, \dots, z_n$ and odd coordinate functions $\zeta_1, \dots, \zeta_r$ . Then we have a holomorphic action of $G$ on $\B$ given by super fractional linear ({\sc Möbius}) transformations

\[
g \rund{\begin{array}{c} \bz \\ \hline
\bzeta \end{array}} := \rund{\begin{array}{c} \rund{A \bz + \bb} \rund{\bc \bz + d}^{- 1} \\ \hline
E \bzeta \rund{\bc \bz + d}^{- 1}
\end{array}} \, ,
\]

where we split

\[
g := \rund{\begin{array}{c|c}
\begin{array}{c|c}
A & \bb \\ \hline
\bc & d\end{array} & 0 \\ \hline
0 & E
\end{array}} \begin{array}{l}
\} n \\
\leftarrow n + 1 \\
\} r
\end{array} \, .
\]

The stabilizer of $\b0 \hookrightarrow \B$ is

\begin{eqnarray*}
&& K := sS\rund{\rund{U(n) \times U(1)} \times U(r)} \\
&& \, = \schweif{\left.\rund{\begin{array}{c|c}
\begin{array}{c|c}
A & 0 \\ \hline
0 & d
\end{array} & 0 \\ \hline
0 & E
\end{array}} \in U(n) \times U(1) \times U(r) \, \right| \, d \det A = \det E}  \, .
\end{eqnarray*}

On $G \times B$ we define the cocycle $j \in \C^\infty(G)^\cz \hat\otimes \O(B)$ as $j(g, \bz) := \rund{\bc \bz + d}^{- 1}$ for all $g \in G$ and $\bz \in B$ . Observe that $j(w) := j(w, \bz) \in U(1)$ is independent of $\bz \in B$ for all $w \in K$ and therefore 
defines a character on the group $K$ . \\

Let $k \in \zz$ be fixed. Then we have a right-representation of $G$

\[
|_g : \D(\B) \rightarrow \D(\B) \, , \, f \mapsto f|_g := f\rund{g \rund{\begin{array}{c}
\bz \\ \hline
\bzeta
\end{array}} } j(g, \bz)^k  \, ,
\]

for all $g \in G$ , which fixes $\O(\B)$ . Finally let $\Gamma$ be a discrete subgroup of~$G$~.

\begin{defin}[super automorphic forms]
Let $f \in \O(\B)$ . Then $f$ is called a super automorphic form for $\Gamma$ of weight $k$ if and only if $f|_\gamma = f$ for all $\gamma \in \Gamma$~. We denote the space of super automorphic forms for $\Gamma$ of weight $k$ by $sM_k(\Gamma)$ .
\end{defin}

Let us define a lift:

\begin{eqnarray*}
\widetilde{\phantom{1}} : \D(\B) &\rightarrow& \C^\infty(G)^\cz \otimes \D\rund{\cz^{0| r}} \simeq \C^\infty(G)^\cz \otimes \bigwedge\rund{\cz^r} \boxtimes \bigwedge\rund{\cz^r} \, , \\
f &\mapsto& \widetilde f  \, ,
\end{eqnarray*}

where

\begin{eqnarray*}
\widetilde f(g) &:=& f|_g\rund{\begin{array}{c}
\b0 \\ \hline
\bEta
\end{array}} \\
&=& f\rund{g \rund{\begin{array}{c}
\b0 \\ \hline
\bEta
\end{array}} } j\rund{g, \b0}^k
\end{eqnarray*}

for all $f \in \D(\B)$ and $g \in G$ and we use the odd coordinate functions $\eta_1, \dots, \eta_r$ on $\cz^{0| r}$ . Let $f \in \O(\B)$ . Then clearly $\widetilde f \in \C^\infty(G)^\cz \otimes \O\rund{\cz^{0| r}}$ and 
$f \in sM_k(\Gamma) \Leftrightarrow \widetilde f \in \C^\infty\rund{\Gamma \backslash G}^\cz \otimes \O\rund{\cz^{0| r}}$ since for all $g \in G$

\[
\begin{array}{ccc}
\C^\infty(G)^\cz \otimes \D\rund{\cz^{0| r}} & \mathop{\longrightarrow}\limits^{l_g} & \C^\infty(G)^\cz \otimes \D\rund{\cz^{0| r}} \\
\uparrow_{\, \widetilde{\phantom{1}}} & \circlearrowleft & \uparrow_{\, \widetilde{\phantom{1}}} \\
\D(\B) & \mathop{\longrightarrow}\limits_{\phantom{1} |_g} & \D(\B)
\end{array} \, ,
\]

where $l_g : \C^\infty(G) \rightarrow \C^\infty(G)$ denotes the left translation with $g \in G$ , $l_g(f)(x) := f(g x)$ for all $x \in G$ . Let $\spitz{\phantom{1}, \phantom{1}}$ be the canonical scalar product on 
$\D\rund{\cz^{0| r}} \simeq \bigwedge\rund{\cz^{2 r}}$ (semilinear in the second entry) . Then for all $a \in \D\rund{\cz^{0| r}}$ we write $\abs{a} := \sqrt{\spitz{a, a}}$ , and $\spitz{\phantom{1}, \phantom{1}}$ induces a 'scalar product'

\[
(f, h)_{\Gamma} := \int_{\Gamma \backslash G} \spitz{\widetilde h, \widetilde f}
\]

for all $f, h \in \D(\B)$ such that $\spitz{\widetilde h, \widetilde f} \in L^1(\Gamma \backslash G)$ and for all $s \in \, ] \, 0, \infty \, ]$ a 'norm'

\[
\norm{f}_{s, \Gamma}^{(k)} := \norm{\phantom{\frac{1}{1}} \abs{\widetilde f} \phantom{\frac{1}{1}}}_{s, \Gamma \backslash G}
\]

for all $f \in \D(\B)$ such that $\abs{\widetilde f} \in \C^\infty\rund{\Gamma \backslash G}$ . On $G$ we always use the (left and right) {\sc Haar} measure. Let us define

\[
L_k^s(\Gamma \backslash \B) := \schweif{f \in \D(\B) \, \left| \begin{array}{c} \\ \\ \end{array}
\widetilde f \in \C^\infty(\Gamma \backslash G)^\cz \otimes \D\rund{\cz^{0|r}} \, , \, \norm{f}_{s, \Gamma}^{(k)} < \infty\right.}   \, .
\]

\begin{defin}[super cusp forms]
Let $f \in sM_k(\Gamma)$ . $f$ is called a super cusp form for $\Gamma$ of weight $k$ if and only if $f \in L_k^2(\Gamma \backslash \B)$ . The
$\cz$- vector space of all super cusp forms for $\Gamma$ of weight $k$ is denoted by $sS_k(\Gamma)$ . It is a {\sc Hilbert} space with inner product $\rund{\phantom{1}, \phantom{1}}_\Gamma$ .
\end{defin}

Observe that $|_g$ respects the splitting

\[
\O(\B) = \bigoplus_{\rho = 0}^r \O^{(\rho)}(\B)
\]

for all $g \in G$ , where $\O^{(\rho)}(\B)$ is the space of all $f = \sum_{I \in \wp(r) \, , \, \abs{I} = \rho} f_I$ , all $f_I \in \O(\B)$ , $I \in \wp(r)$ , $\abs{I} = \rho$ , $\rho = 0, \dots, r$ , and $\, \widetilde{\phantom{1}} \,$ maps the space $\O^{(\rho)}(\B)$ into
$\C^\infty(G)^\cz \otimes \O^{(\rho)}\rund{\cz^{0| r}}$ . Therefore we have splittings

\[
sM_k(\Gamma) = \bigoplus_{\rho = 0}^r sM_k^{(\rho)}(\Gamma) \phantom{1} \text{ and } \phantom{1} sS_k(\Gamma) = \bigoplus_{\rho = 0}^r sS_k^{(\rho)}(\Gamma)  \, ,
\]

where $sM_k^{(\rho)}(\Gamma) := sM_k(\Gamma) \cap \O^{(\rho)}(\B)$ , $sS_k^{(\rho)}(\Gamma) := sS_k(\Gamma) \cap \O^{(\rho)}(\B)$ , $\rho = 0, \dots, r$ , and the last sum is orthogonal. \\

As I show in \cite{Kne} and in section 3.2 of \cite{KneBuch} there is an analogon to {\sc Satake}'s theorem in the super case:

\begin{theorem} \label{Satake} Let $\rho \in \{0, \dots, r\}$ . Assume $\Gamma \backslash G$ is compact or $n \geq 2$ and $\Gamma \sqsubset G$ is a lattice (discrete such that $\vol \, \Gamma \backslash G < \infty$ , $\Gamma \backslash G$ not 
necessarily compact) . If $k \geq 2 n - \rho$ then

\[
sS_k^{(\rho)}(\Gamma) = sM_k^{(\rho)}(\Gamma) \cap L_k^s \rund{\Gamma \backslash \B}
\]

for all $s \in \, [ \, 1, \infty \, ] \,$ .
\end{theorem}

As in the classical case this theorem implies that if $\Gamma \backslash G$ is compact or $n \geq 2$ , $\Gamma \sqsubset G$ is a lattice and $k \geq 2 n - \rho$ then the {\sc Hilbert} space $sS_k^{(\rho)}(\Gamma)$ is finite dimensional. \\

We will use the {\sc Jordan} triple determinant $\Delta: \cz^n \times \cz^n \rightarrow \cz$ given by

\[
\Delta\rund{\bz, \bw} := 1 - \bw^* \bz
\]

for all $\bz, \bw \in \cz^n$ . Let us recall the basic properties:

\begin{itemize}
\item[(i)] $\abs{j\rund{g, \b0}} = \Delta\rund{g \b0, g \b0}^\frac{1}{2}$ for all $g \in G$ ,

\item[(ii)] $\Delta\rund{g \bz, g \bw} = \Delta\rund{\bz, \bw} j\rund{g, \bz} \overline{j\rund{g, \bw}}$ for all $g \in G$ and $\bz, \bw \in B$ , and

\item[(iii)] $\int_B \Delta\rund{\bz, \bz}^\lambda d V_{\Leb} < \infty$ if and only if $\lambda > - 1$ .
\end{itemize}

We have the $G$-invariant volume element $\Delta(\bz, \bz)^{- (n + 1)} d V_{\Leb}$ on $B$ . \\

For all $I \in \wp(r)$ , $h \in \O(B)$ , $\bz \in B$ and \\
$g = \rund{\begin{array}{c|c}
* & 0 \\ \hline
0 & E
\end{array}} \in G$ we have

\[
\left.h \bzeta^I \right|_g \rund{\bz} = h\rund{g \bz} \rund{E \bEta}^I j\rund{g, \bz}^{k + \abs{I}} \, ,
\]

where $E \in U(r)$ . So for all $s \in \, ] \, 0, \infty \, ] \, $ , $f = \sum_{I \in \wp(r)} f_I \bzeta^I$ and \\
$h = \sum_{I \in \wp(r)} h_I \bzeta^I \in \O(\B)$ we have

\[
\norm{f}_{s, \Gamma}^{(k)} \equiv \norm{\, \sqrt{\sum_{I \in \wp(r)} f_I^2 \Delta\rund{\bz, \bz}^{k + \abs{I}} } \, }_{s, \,  \Gamma \backslash B, \, \Delta\rund{\bz, \bz}^{- (n + 1)} d V_{\Leb}}
\]

if $\widetilde f \in \C^\infty(G) \otimes \O\rund{\cz^{0| r}}$ and

\[
(f, h)_\Gamma \equiv \sum_{I \in \wp(r)} \int_{\Gamma \backslash B} \overline{f_I} h_I \Delta\rund{\bz, \bz}^{k + \abs{I} - (n + 1)} d V_{\Leb}
\]

if $\spitz{\widetilde h , \widetilde f} \in L^1(\Gamma \backslash G)$ , where '$\equiv$' means equality up to a constant $\not= 0$ depending on $\Gamma$ . \\

For the explicit computation of the elements of our spanning set in theorem \ref{main} we need the following lemmas:

\begin{lemma}[convergence of relative {\sc Poincaré} series] \label{Poincare} Let $\Gamma_0 \sqsubset \Gamma$ be a subgroup and

\[
f \in sM_k\rund{\Gamma_0} \cap L_k^1 \rund{\Gamma_0 \backslash \B} \, .
\]

Then

\[
\Phi := \sum_{\gamma \in \Gamma_0 \backslash \Gamma} f |_\gamma \text{  and  }  \Phi' := \sum_{\gamma \in \Gamma_0 \backslash \Gamma} \widetilde f(\gamma \diamondsuit)
\]

converge absolutely and uniformly on compact subsets of $B$ resp. $G$ ,

\[
\Phi \in sM_k (\Gamma) \cap L_k^1 \rund{\Gamma \backslash \B} \,  ,
\]

$\widetilde \Phi = \Phi'$ , and for all $\varphi \in sM_k (\Gamma) \cap L_k^\infty \rund{\Gamma \backslash \B}$ we have

\[
\rund{\Phi, \varphi}_\Gamma = \rund{f, \varphi}_{\Gamma_0} \,  .
\]

\end{lemma}

The symbol '$\diamondsuit$' here and also later simply stands for the argument of the function. So $\widetilde f (\gamma \diamondsuit) \in \C^\infty(G)^\cz \otimes \bigwedge\rund{\cz^r}$ is a short notation for the smooth map

\[
G \rightarrow \bigwedge\rund{\cz^r} \, , \, g \mapsto \widetilde f(\gamma g) \, .
\]

{\it Proof:} standard using the mean value property of holomorphic functions for all $k \in \zz$ without any further assumption on $k$ . $\Box$ \\

\begin{lemma} \label{reproducing}
Let $I \in \wp(r)$ and $k \geq 2 n + 1 - \abs{I}$ . Then for all $\bw \in B$

\[
\Delta\rund{\diamondsuit, \bw}^{- k - \abs{I}} \bzeta^I \in \O^{\abs{I}}(\B) \cap L_k^1(\B)  \, ,
\]

and for all $f = \sum_{J \in \wp(r)} f_J \bzeta^J \in \O(\B) \cap L_k^\infty(\B)$ we have

\[
\rund{\Delta\rund{\diamondsuit, \bw}^{- k - \abs{I}} \bzeta^I, f} \equiv f_I\rund{\bw} \, ,
\]

where $\rund{\phantom{1}, \phantom{1}} := \rund{\phantom{1}, \phantom{1}}_{\{1\}}$ .

\end{lemma}

Since the proof is also standard, we will omit it here. It can be found in~\cite{KneBuch}~.

\section{The structure of the group $G$ }

We have a canonical embedding

\[
G' := SU(p, q) \hookrightarrow G \, , \, g' \mapsto \rund{\begin{array}{c|c}
g' & 0 \\ \hline
0 & 1
\end{array}} \,  ,
\]

and the canonical projection

\[
G \rightarrow U(r) \, , \, g := \rund{\begin{array}{c|c}
g' & 0 \\ \hline
0 & E
\end{array}} \mapsto E_g := E
\]

induces a group isomorphism

\[
G \left/ G' \right. \simeq U(r)  \, .
\]

Obviously $K' = K \cap G' = S(U(n) \times U(1))$ is the stabilizer of $\b0$ in $G'$ . Let $A$ denote the common standard maximal split abelian subgroup of $G$ and $G'$ given by the image of the {\sc Lie} group embedding

\[
\rz \hookrightarrow G' \, , \, t \mapsto a_t :=
\rund{\begin{array}{c|c}

\begin{array}{c|c}
\cosh t & 0 \\ \hline
0 & 1
\end{array} & \begin{array}{c}
\sinh t_1 \\ \hline
0
\end{array} \\ \hline

\begin{array}{c|c} \sinh t & 0 \end{array} & \cosh t

\end{array} }  \,  .
\]

Then the centralizer $M$ of $A$ in $K$ is the group of all

\[
\rund{\begin{array}{c|c}
\begin{array}{c|c}
\begin{array}{c|c}
\eps & 0 \\ \hline

0 & u \end{array} & 0 \\ \hline
0 & \eps \end{array} & 0 \\ \hline
0 & E
\end{array}} \,  ,
\]

where $\eps \in U\rund{1}$ , $u \in U\rund{p - 1}$ and $E \in U(r)$ such that $\eps^2 \det u = \det E$ . Let $M' = K' \cap M = G' \cap M$ be the centralizer of $A$ in $K'$ . The centralizer of $G'$ in $G$ is precisely

\[
Z_G\rund{G'} := \schweif{\left.\rund{\begin{array}{c|c}
\eps 1 & 0 \\ \hline
0 & E
\end{array}} \, \right| \, \eps \in U(1) \, , \, E \in U(r) \, , \, \eps^{p + 1} = \det E} \sqsubset M  \, ,
\]

and $G' \cap Z_G\rund{G'} = Z\rund{G'}$ . An easy calculation shows that $G = G' Z_G\rund{G'}$~. So $K = K' Z_G\rund{G'}$ and $M = M' Z\rund{G'}$ . Therefore if we decompose the adjoint representation of $A$ as

\[
\g = \bigoplus_{\alpha \in \Phi} \g^{\alpha} \,  ,
\]

where for all $\alpha \in \rz$

\[
\g^{\alpha} := \schweif{\xi \in \g \, \left| \, \Ad_{a_t} (\xi) = e^{\alpha t} \right. }
\]

is the corresponding root space and

\[
\Phi := \schweif{\left. \alpha \in \rz \, \right| \, g^{\alpha} \not= 0}
\]

is the root system, then we see that $\Phi$ is at the same time the root system of $G'$ , so $\Phi = \schweif{0, \pm 2}$ if $n = 1$ and $\Phi = \schweif{0, \pm 1, \pm 2}$ if $n \geq 2$ , furthermore if $\alpha \not= 0$ then $\g^\alpha \sqsubset \g'$ is at the
same time the corresponding root space of $\g'$~, and finally $\g^0 = \a \oplus \m = \a \oplus \m' \oplus \z_{\g}\rund{\g'}$ .

\begin{lemma} \label{normalizer}

\[
N(A) = A N_K(A) = N(A M) \sqsubset N(M) \,  .
\]

\end{lemma}

{\it Proof:} simple calculation. $\Box$ \\

In particular we have the {\sc Weyl} group

\[
W := M \left\backslash N_K(A) \right. \simeq M' \left\backslash N_{K'}(A) \right. \simeq \{\pm 1\}
\]

acting on $A \simeq \rz$ via sign change. For the main result, theorem \ref{main} , of this article the following definition is crucial:

\begin{defin}
Let $g_0 \in G$ .

\item[(i)] $g_0$ is called loxodromic if and only if there exists $g \in G$ such that\\
$g_0 \in g A M g^{- 1}$ .
\item[(ii)] If $g_0$ is loxodromic, it is called regular if and only if $g_0 = g a_t w g^{- 1}$ with $t \in \rz \setminus \{0\}$ and $w \in M$ .
\item[(iii)] If $\gamma \in \Gamma$ is regular loxodromic then it is called primitive in $\Gamma$ if and only if $\gamma = \gamma'^\nu$ implies $\nu \in \{\pm 1\}$ for all loxodromic $\gamma' \in \Gamma$ and $\nu \in \zz$ .
\end{defin}

Clearly for all $\gamma \in \Gamma$ regular loxodromic there exists $\gamma' \in \Gamma$ primitive regular loxodromic and $\nu \in \nz \setminus \{0\}$ such that $\gamma = \gamma'^\nu$ .

\begin{lemma} \label{determined} Let $g_0 \in G$ be regular loxodromic, $g \in G$ , $w \in M$ and \\
$t \in \rz \setminus~\{0\}$ such that $g_0 = g a_t w g^{- 1}$ . Then $g$ is uniquely determined up to right translation by elements of $A N_K(A)$ , and $t$ is uniquely determined up to sign.
\end{lemma}

{\it Proof:} by straight forward computation or using the following trick: Let $g' \in G$ , $w' \in M$ and $t' \in \rz$ such that also $g_0 = g' a_{t'} w' g'^{- 1}$ . Then $a_t w = \rund{g^{- 1} g'} a_{t'} w' \rund{g^{- 1} g'}^{- 1}$ . Since $t \in \rz \setminus \{0\}$ and
because of the root space decomposition, $\a + \m$ must be the largest subspace of $\g$ on which $\Ad_{a_t w}$ is orthogonal with respect to an appropiate scalar product. So $\Ad_{g^{- 1} g'}$ maps $\a + \m$ into itself. This implies
$g^{- 1} g' \in N (A M) = A N_K(A)$ by lemma \ref{normalizer}~.~$\Box$

\section{The main result}

Let $\rho \in \{0, \dots, r\}$ . Assume $\Gamma \backslash G$ compact {\bf or} $n \geq 2$ , $\vol \Gamma \backslash G < \infty$ and $k \geq 2 n - \rho$ . Let $C > 0$ be given. Let us consider a regular loxodromic $\gamma_0 \in \Gamma$ . Let $g \in G$ ,
$w_0 \in M$ and $t_0 > 0$ such that $\gamma_0 = g a_{t_0} w_0 g^{- 1}$ . \\

There exists a torus $\tz := \left.\spitz{\gamma_0} \right\backslash g A M$ belonging to $\gamma_0$ . From lemma \ref{determined} it follows that $\tz$ is independent of $g$ up to right translation with an element of the {\sc Weyl} group 
$W = M \backslash N_K(A)$ . \\

Let $f \in sS_k(\Gamma)$ . Then $\widetilde f \in \C^\infty \rund{\Gamma \backslash G}^\cz \otimes \O\rund{\cz^{0| r}}$ . Define \\
$h \in \C^\infty \rund{\rz \times M}^\cz \otimes \O\rund{\cz^{0| r}}$ as

\[
h \rund{t, w} := \widetilde f \rund{g a_t w}
\]

for all $(t, w) \in \rz \times M$ 'screening' the values of $\widetilde f$ on $\tz$ . Then clearly $h \rund{t, w} = h\rund{t, 1, E_w \bEta j(w)} j(w)^k$ , and so $h(t, w) = h(t, 1, E_w \eta) j(w)^{k + \rho}$ if $f \in sS_k^{(\rho)}(\Gamma)$ , for all
$\rund{t, w} \in \rz \times M$ . Clearly $E_0 := E_{w_0} \in U(r)$ . So we can choose $g \in G$ such that $E_0$ is diagonal without changing $\tz$ . Choose $D \in \rz^{r \times r}$ diagonal such that $\exp(2 \pi i D) = E_0$
and $\chi \in \rz$ such that $j\rund{w_0} = e^{2 \pi i \chi}$ . $D$ and $\chi$ are uniquely determined by $w_0$ up to $\zz$ . If $D = \rund{\begin{array}{ccc}
d_1 & & 0 \\
 & \ddots & \\
0 & & d_r
\end{array}}$ with $d_1 , \dots, d_r \in \rz$ and $I \in \wp(r)$ then we define $\tr_I D := \sum_{j \in I} d_j$ . \\

\begin{theorem}[{\sc Fourier} expansion of $h$ ] \label{Fourier}

\item[(i)] $h \rund{t + t_0 , w} = h \rund{t , w_0^{- 1} w}$ for all $\rund{t, w} \in \rz \times M$ , and there exist unique $b_{I, m} \in \cz$ , $I \in \wp(r)$ ,
$m \in \frac{1}{t_0} \rund{\zz - \rund{k + \abs{I}} \chi - \tr_I D}$ , such that

\[
h \rund{t, w} = \sum_{I \in \wp(r)} j(w)^{k + \abs{I}} \sum_{m \in \frac{1}{t_0} \rund{\zz - \rund{k + \abs{I}} \chi - \tr_I D}} b_{I, m} e^{2 \pi i m t} \rund{E_w \bEta}^I
\]

for all $\rund{t, w} \in \rz \times M$ , where the sum converges uniformly in all derivatives.

\item[(ii)] If $f \in sS_k^{(\rho)}(\Gamma)$ , $b_{I, m} = 0$ for all $I \in \wp(r)$ , $\abs{I} = \rho$ , and \\
$m \in \frac{1}{t_0} \rund{\zz - \rund{k + \rho} \chi - \tr_I D} \cap \, ] \, - C, C \, [$ then there exists \\
$H \in \C^\infty \rund{\rz \times M}^\cz \otimes \bigwedge\rund{\cz^r}$ uniformly {\sc Lipschitz} continuous with a {\sc Lipschitz} constant $C_2 \geq 0$ independent of $\gamma_0$ such that

\[
h = \partial_t H \, ,
\]

\[
H\rund{t, w} = j(w)^k H\rund{t, 1, E_w \bEta j(w)}
\]

and

\[
H\rund{t + t_0 , w} = H\rund{t , w_0^{- 1} w}
\]

for all $\rund{t, w} \in \rz \times M$ .

\end{theorem}

{\it Proof:} (i) Let $t \in \rz$ and $w \in M$ . Then

\begin{eqnarray*}
h\rund{t + t_0 , w} &=& \widetilde f \rund{g a_{t_0} a_t w} = \widetilde f \rund{\gamma_0 g w_0^{- 1} a_t w} = \widetilde f \rund{g a_t w_0^{- 1} w} \\
&=& h\rund{t, w_0^{- 1} w}  \, ,
\end{eqnarray*}

and so

\begin{eqnarray*}
h\rund{t + t_0 , 1} &=& h\rund{t, w_0^{- 1}} \\
&=& j\rund{w_0}^{- k} h\rund{t, 1, E_0^{- 1} \bEta j\rund{w_0}^{- 1}} \\
&=& j\rund{w_0}^{- k} \sum_{I \in \wp(r)} h\rund{t, 1} e^{- 2 \pi i \tr_I D} \bEta^I j\rund{w_0}^{- \abs{I}} \\
&=& \sum_{I \in \wp(r)} e^{- 2 \pi i \rund{\rund{k + \abs{I}} \chi + \tr_I D}} h_I(t, 1) \bEta^I  \, .
\end{eqnarray*}

Therefore $h_I\rund{t + t_0 , 1} = e^{- 2 \pi i \rund{\rund{k + \abs{I}} \chi + \tr_I D}} h_I(t, 1)$ for all $I \in \wp(r)$ , and the rest follows by standard {\sc Fourier} expansion. $\Box$ \\

For proving (ii) we need the following lemma:

\begin{lemma}[generalization of the reverse {\sc Bernstein} inequality] \label{Bernstein}

Let $t_0 \in \rz \setminus \{0\}$ , $\nu \in \rz$ and $C > 0$ . Let $\S$ be the space of all convergent {\sc Fourier} series

\[
s = \sum_{m \in \frac{1}{t_0} (\zz - \nu) \, , \, \abs{m} \geq C} s_l e^{2 \pi i m \diamondsuit} \in \C^\infty\rund{\rz}^\cz  \, ,
\]

all $s_m \in \cz$ . Then

\[
\widehat{\phantom{1}} : \S \rightarrow \S \, , \, s = \sum_{m \in \frac{1}{t_0} (\zz - \nu) \, , \, \abs{m} \geq C} s_m e^{2 \pi i m \diamondsuit} \mapsto
\widehat s := \sum_{m \in \frac{1}{t_0} (\zz - \nu) \, , \, \abs{m} \geq C} \frac{s_m}{2 \pi i m} e^{2 \pi i m \diamondsuit}
\]

is a well-defined linear map, and $\norm{\widehat s}_\infty \leq \frac{6}{\pi C} \norm{s}_\infty$ for all $s \in \S$ .
\end{lemma}

{\it Proof:} This can be deduced from the ordinary reverse {\sc Bernstein} inequality, see for example theorem 8.4 in
chapter I of \cite{Katz} . $\Box$ \\

{\it Now we prove theorem \ref{Fourier} (ii) .} Fix some $I \in \wp(r)$ such that $\abs{I} = \rho$ and $b_{I, m} = 0$ for all $m \in \frac{1}{t_0} \rund{\zz - \rund{k + \rho} \chi - \tr_I D} \cap \, ] \, - C, C \, [$ . Then if we define 
$\nu := \rund{k + \rho} \chi + \tr_I D \in \rz$ we have

\[
h_I(\diamondsuit, 1) = \sum_{m \in \frac{1}{t_0} \rund{\zz - \nu} \, , \, \abs{m} \geq C} b_{I, m} e^{2 \pi i m \diamondsuit}  \, ,
\]

and so we can apply the generalized reverse {\sc Bernstein} inequality, lemma \ref{Bernstein} , to $h_I$ . Therefore we can define

\[
H_I' := \widehat{h_I \rund{\diamondsuit, 1}} = \sum_{m \in \frac{1}{t_0} \rund{\zz - \nu} \, , \, \abs{m} \geq C} \frac{b_{I, m}}{2 \pi i m} e^{2 \pi i m \diamondsuit}
\in \C^\infty \rund{\rz}^\cz \, .
\]

$\abs{\widetilde f} \in L^\infty (G)$ by {\sc Satake}'s theorem , theorem \ref{Satake} , and so there exists a constant $C' > 0$ independent of $\gamma_0$ and $I$ such that $\norm{h_I}_\infty < C'$ , and now lemma \ref{Bernstein} tells us that

\[
\norm{H_I'}_\infty \leq \frac{6}{\pi C} \norm{h \rund{\diamondsuit, 1}}_\infty \leq \frac{6 C'}{\pi C}  \, .
\]

Clearly $h_I \rund{\diamondsuit, 1} = \partial_t H_I'$ . \\

Since $j$ is smooth on the compact set $M$ , $j^{k + \rho} \rund{E_w \bEta}^I$ is uniformly {\sc Lipschitz} continuous on $M$ with a common {\sc Lipschitz} constant $C''$ independent of $\gamma_0$ and $I$ . So we see that 
$H \in \C^\infty(\rz, M)^\cz \otimes \bigwedge\rund{\cz^r}$ defined as

\[
H(t, w) := \sum_{I \in \wp(r)} j(w)^{k + \rho} H_I'(t)\rund{E_w \bEta}^I
\]

for all $(t, w) \in \rz \times M$ is uniformly {\sc Lipschitz} continuous with {\sc Lipschitz} constant $C_2 := \rund{\frac{6 C''}{\pi C} + 1} C'$ independent of $\gamma_0$ , and the rest is trivial.
$\Box$ \\

Let $I \in \wp(r)$ and $m \in \frac{1}{t_0} \rund{\zz - \rund{k + \abs{I}} \chi - \tr_I D}$ . Since $sS_k (\Gamma)$ is a {\sc Hilbert} space and
$sS_k (\Gamma) \rightarrow \cz \, , \, f \mapsto b_{I, m}$ is linear and continuous there exists exactly one $\varphi_{\gamma_0 , I, m} \in sS_k (\Gamma)$ such that
$b_{I, m} = \rund{\varphi_{\gamma_0 , I, m}, f}$ for all $f \in sS_k(\Gamma)$ . Clearly $\varphi_{\gamma_0 , I, m} \in sS_k^{\rund{\abs{I}} } (\Gamma)$ . \\

From now on for the rest of the article for simplicity we write $m \in \, ] \, - C, C \, [$ instead of $m \in \frac{1}{t_0} \rund{\zz - \rund{k + \abs{I}} \chi - \tr_I D} \cap \, ] \, - C, C \, [$ . In the last section we will compute $\varphi_{\gamma_0 , I, m}$ as a relative
{\sc Poincaré} series. One can check that the family

\[
\schweif{\varphi_{\gamma_0 , I, m}}_{I \in \wp(r) \, , \abs{I} = \rho \, , \, m \in  \, ] \, - C, C \, [}
\]

is independent of the choice of $g$ , $D$ and $\chi$ up to multiplication with a unitary matrix with entries in $\cz$ and invariant under conjugating $\gamma_0$ with elements of $\Gamma$ . \\

Now we can state our main theorem: Let $\Omega$ be a fundamental set for all primitive regular loxodromic $\gamma_0 \in \Gamma$ modulo conjugation by elements of $\Gamma$  and

\[
\widetilde Z := \overline{\schweif{\left.m \in Z_G\rund{G'} \phantom{\frac{}{}} \right| \phantom{\frac{\frac{}{}}{}} \exists \, g' \in G' \, : \, m g' \in \Gamma}} \sqsubset Z_G\rund{G'}  \, .
\]

Then clearly $\Gamma \sqsubset G' \widetilde Z$ . {\bf Recall that we still assume}

\begin{itemize}
\item $\Gamma \backslash G$ {\bf compact} or \\
\item $n \geq 2$ , $\vol \Gamma \backslash G < \infty$ and $k \geq 2 n - \rho$ .
\end{itemize}

\begin{theorem}[spanning set for $sS_k (\Gamma)$ ] \label{main} Assume that the right translation of $A$ on $\Gamma \backslash G' \widetilde Z$ is topologically transitive. Then

\[
\schweif{\left. \varphi_{\gamma_0 , I, m} \, \right| \, \gamma_0 \in \Omega , I \in \wp(r) , \abs{I} = \rho \, , \, m \in \, ] \, - C, C \, [ \, }
\]

is a spanning set for $sS_k^{(\rho)} (\Gamma)$ .
\end{theorem}

For proving this result we need an {\sc Anosov} type theorem for $G$ and the unbounded realization of $\B$ , which we will discuss in the following two sections. \\

{\it Remarks:}

\begin{itemize}
\item[(i)] If there is some subgroup $\widetilde M \sqsubset Z_G\rund{G'}$ such that $\Gamma \sqsubset G' \widetilde M$ and the right translation of $A$ on $\Gamma \backslash G' \widetilde M$ is topologically transitive then necessarily 
$\widetilde M Z(G') = \widetilde Z$ and there exists $g_0 \in G'$ such that \\
$G' \widetilde Z = \overline{\Gamma g_0 A}$ . The latter statement is a trivial consequence of the fact that $\widetilde Z \sqsubset M$ .
\item[(ii)] In the case where $\Gamma \cap G' \sqsubset \Gamma$ is of finite index or equivalently $\widetilde Z$ is finite then we know that the right translation of $A$ on $\Gamma \backslash G' \widetilde Z$ is topologically transitive because of 
{\sc Moore}'s ergodicity theorem, see \cite{Zim} theorem 2.2.6 , and since then $\Gamma \cap G' \sqsubset G'$ is a lattice.
\item[(iii)] There is a finite-to-one correspondence between $\Omega$ and the set of closed geodesics of $\Gamma \backslash B$ assigning to each primitive loxodromic element \\
$\gamma_0 = g a_{t_0} w_0 g^{- 1} \in \Gamma$ , $g \in G$ , $t_0 > 0$ and $w_0 \in M$ , the image of the unique geodesic $g A \b0$ of $B$ normalized by $\gamma_0$ under the canonical projection $B \rightarrow \Gamma \backslash B$ . It is of length 
$t_0$ if there is no irregular point of $\Gamma \backslash B$ on $g A \b0$ .
\end{itemize}

\section{An {\sc Anosov} type result for the group $G$ }

On the {\sc Lie} group $G$ we have a smooth flow $\rund{\varphi_t}_{t \in \rz}$ given by the right translation by elements of $A$ :

\[
\varphi_t : G \rightarrow G \, , \, g \mapsto g a_t \,  .
\]

This turns out to be partially hyperbolic, and so we can apply a partial {\sc Anosov} closing lemma. By the way the flow $\rund{\varphi_t}_{t \in \rz}$ descends to the ordinary geodesic flow on the unit tangent bundle $S B \simeq G / M$ . Let us 
first have a look at the general theory of partial hyperbolicity: Let $W$ be for the moment a smooth Riemannian manifold.

\begin{defin}[partially hyperbolic diffeomorphism and flow] Let $C > 1$ .
\item[(i)] Let $\varphi$ be a $\C^\infty$-diffeomorphism of $W$ . Then $\varphi$ is called partially hyperbolic with constant $C$ if and only if there exists an orthogonal $D \varphi$ (and therefore $D \varphi^{- 1}$ ) -invariant $\C^\infty$-splitting

\begin{equation}
T W = T^0 \oplus T^+ \oplus T^-  \label{splitting}
\end{equation}

of the tangent bundle $T W$ such that $T^0 \oplus T^+$ , $T^0 \oplus T^-$ , $T^0$ , $T^+$ and $T^-$ are closed under the commutator, $D \varphi |_{T^0}$ is an isometry,
$\norm{D \varphi |_{T^-}} \leq \frac{1}{C}$ and $\norm{D \varphi^{- 1} |_{T^+}} \leq \frac{1}{C}$ .
\item[(ii)] Let $\rund{\varphi_t}_{t \in \rz}$ be a $\C^\infty$-flow on $W$ . Then $\rund{\varphi_t}_{t \in \rz}$ is called partially hyperbolic with constant $C$ if and only if all $\varphi_t$ , $t > 0$ are partially hyperbolic diffeomorphisms with a common splitting
(\ref{splitting}) and constants $e^{C t}$ resp. and $T^0$ contains the generator of the flow.
\end{defin}

A partially hyperbolic diffeomorphism $\varphi$ gives rise to $\C^\infty$-foliations on $W$ corresponding to the splitting $T W = T^0 \oplus T^+ \oplus T^-$ . Let us denote the distances along the
$T^0 \oplus T^+$- , $T^0$-~, $T^+$- respectively $T^-$-leaves by $d^{0, +}$ , $d^0$~, $d^+$ and $d^-$ .

\begin{defin} \label{SU(p,q)hyper rectangular}
Let $T W = T^0 \oplus T^+ \oplus T^-$ be an orthogonal $\C^\infty$-splitting of the tangent bundle $T W$ of $W$ such that $T^0 \oplus T^+$ , $T^0$ , $T^+$ and $T^-$ are closed
under the commutator, $C' \geq 1$ and $U \subset W$ . $U$ is called $C'$-rectangular (with respect to the splitting $T W = T^0 \oplus T^+ \oplus T^-$ ) if and only if for all $y , z \in U$

\begin{itemize}
\item[\{i\}] there exists a unique intersection point $a \in U$ of the $T^0 \oplus T^+$-leaf containing $y$ and the $T^-$-leaf containing $z$ and a unique intersection point $b \in U$ of the
$T^0 \oplus T^+$-leaf containing $z$ and the $T^-$-leaf containing~$y$~,

\[
d^{0,+}\rund{y, a} , d^-\rund{y, b} , d^-\rund{z, a} , d^{0, +}\rund{z, b} \leq C' d\rund{y, z} \, ,
\]

and

\begin{eqnarray*}
\frac{1}{C'} d^{0, +}\rund{z, b} \leq d^{0, +}\rund{y, a} \leq C' d^{0, +}\rund{z, b} , \\
\frac{1}{C'} d^-\rund{z, a} \leq d^-\rund{y, b} \leq C' d^-\rund{z, a} \, .
\end{eqnarray*}

\item[\{ii\}]

if $y$ and $z$ belong to the same $T^0 \oplus T^+$-leaf there exists a unique intersection point $c \in U$ of the $T^0$-leaf containing $y$ and the $T^+$-leaf containing $z$ and a unique
intersection point $d \in U$ of the $T^0$-leaf containing $z$ and the $T^+$-leaf containing $y$ ,

\[
d^0\rund{y, c} , d^+\rund{y, d} , d^+\rund{z, c} , d^0\rund{z, d} \leq C' d^{0, +}\rund{y, z} \, ,
\]

and

\begin{eqnarray*}
\frac{1}{C'} d^0\rund{z, d} \leq d^0\rund{y, c} \leq C' d^0\rund{z, d} \, , \\
\frac{1}{C'} d^+\rund{z, c} \leq d^+\rund{y, d} \leq C' d^+\rund{z, c} \, .
\end{eqnarray*}

\end{itemize}

\begin{figure}[H]
\begin{center}
\includegraphics[width=0.8\textwidth]{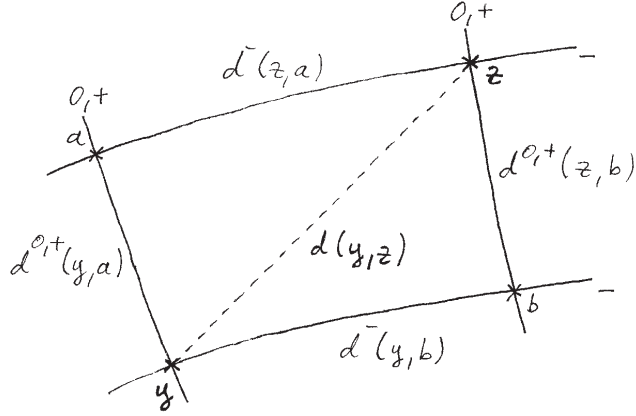}
\caption{intersection points in \{i\} .}
\end{center}
\end{figure}

\end{defin}

Since the splitting $T W = T^0 \oplus T^+ \oplus T^-$ is orthogonal and smooth we see that for all $x \in W$ and $C' > 1$ there exists a $C'$-rectangular neighbourhood of~$x$~.

\begin{theorem}[partial {\sc Anosov} closing lemma] \label{closing}

Let $\varphi$ be a partially hyperbolic diffeomorphism with constant $C$ , let $x \in W$ , $C' \in \, ] \, 1, C \, [$ and $\delta > 0$ such that $\overline{U_\delta(x)}$ is contained in a $C'$-rectangular subset $U \subset W$ .

If $d\rund{x , \varphi(x)} \leq \delta \frac{1 - \frac{C'}{C}}{C'^2 +1}$ then there exist $y, z \in U$ such that

\begin{itemize}

\item[(i)] $x$ and $y$ belong to the same $T^-$-leaf and

\[
d^-\rund{x, y} \leq \frac{C'}{1 - \frac{C'}{C}} d\rund{x , \varphi(x)} \, ,
\]

\item[(ii)] $y$ and $\varphi(y)$ belong to the same $T^0 \oplus T^+$-leaf and

\[
d^{0, +} \rund{y, \varphi(y)} \leq C'^2 d\rund{x , \varphi(x)} \,  ,
\]

\item[(iii)] $y$ and $z$ belong to the same $T^+$-leaf and

\[
d^+\rund{\varphi(y), \varphi(z)} \leq \frac{C'^3}{1 - \frac{C'}{C}} d\rund{x , \varphi(x)} \, ,
\]

\item[(iv)] $z$ and $\varphi(z)$ belong to the same $T^0$-leaf and

\[
d^0\rund{z, \varphi(z)} \leq \C'^4 d\rund{x , \varphi(x)} \, .
\]

\end{itemize}
\end{theorem}

The proof, which will not be given here, uses a standard argument obtaining the points $y$ and $\varphi(z)$ as limits of certain {\sc Cauchy} sequences. The interested reader will find it in \cite{KneBuch} . \\

Now let us return to the flow $\rund{\varphi_t}_{t \in \rz}$ on $G$ and choose a left invariant metric on $G$ such that $\g^\alpha$ , $\alpha \in \Phi \setminus \{0\}$ , $\a$ and $\m$ are pairwise orthogonal and the isomorphism $\rz \simeq A \subset G$ is
isometric. Then since the flow $\rund{\varphi_t}_{t \in \rz}$ commutes with left translations it is indeed partially hyperbolic with constant $1$ and the unique left invariant splitting of $T G$ given by

\[
T_1 G = \g = \underbrace{\begin{array}{c}
\\
\\
\end{array} \a \oplus \m \begin{array}{c}
\\
\\
\end{array}}_{T_1^0 :=} \oplus \underbrace{\bigoplus_{\alpha \in \Phi \, , \, \alpha > 0} \g^{\alpha}}_{T_1^- :=}
\oplus \underbrace{\bigoplus_{\alpha \in \Phi \, , \, \alpha < 0} \g^{\alpha}}_{T_1^+ :=} \,  .
\]

For all $L \subset G$ compact, $T, \eps > 0$ define

\[
M_{L, T} := \schweif{\left. g a_t g^{- 1} \, \right| \, g \in L , t \in [- T, T]}
\]

and

\[
N_{L, T, \eps} := \schweif{g \in G \, \left| \, \dist\rund{g, M_{L, T}} \leq \eps \right.} \,  .
\]

\begin{lemma} \label{Gamma}
For all $L \subset G$ compact there exist $T_0 , \eps_0 > 0$ such that\\
$\Gamma \cap N_{L, T_0, \eps_0} = \{1\}$ .
\end{lemma}

{\it Proof:} Let $L \subset G$ be compact and $T > 0$ . Then $M_{L, T}$ is compact, and so there exists $\eps > 0$ such that $N_{L, T, \eps}$ is again compact. Since $\Gamma$ is
discrete, $\Gamma \cap N_{L, T, \eps}$ is finite. Clearly for all $T, T', \eps$ and $\eps' > 0$ if $T \leq T'$ and $\eps \leq \eps'$ then $N_{L, T, \eps} \subset N_{L, T', \eps'}$ , and finally

\[
\bigcap_{T, \eps > 0} N_{T, \eps} = \{1\} \, . \, \Box
\]

Here now the quintessence of this section:

\pagebreak

\begin{theorem} \label{Anosov}

\item[(i)] For all $T_1 > 0$ there exist $C_1 \geq 1$ and $\eps_1 > 0$ such that for all $x \in G$ , $\gamma \in \Gamma$ and $T \geq T_1$ if

\[
\eps := d\rund{\gamma x, x a_T} \leq \eps_1
\]

then there exist $z \in G$ , $w \in M$ and $t_0 > 0$ such that $\gamma z = z a_{t_0} w$ (and so $\gamma$ is regular loxodromic) , $d\rund{(t_0, w) , (T, 1)} \leq C_1 \eps$ and for all $\tau \in \, [ \, 0, T \, ]$

\[
d\rund{x a_\tau , z a_\tau} \leq C_1 \eps \rund{e^{- \tau} + e^{- \rund{T - \tau}} }  \, .
\]

\item[(ii)] For all $L \subset G$ compact there exists $\eps_2 > 0$ such that for all $x \in L$ , $\gamma \in \Gamma$ and $T \in \, \eckig{ \, 0, T_0 \, }$ , $T_0 > 0$ given by lemma
\ref{Gamma} , if

\[
\eps := d\rund{\gamma x, x a_T} \leq \eps_2
\]

then $\gamma = 1$ and $T \leq \eps$ .

\end{theorem}

{\it Proof:} (i) Let $T_1 > 0$ and define

\[
C_1 := \max\rund{\frac{e^{\frac{3}{2} T_1}}{1 - e^{- \frac{T_1}{2}} } , e^{2 T_1} } \geq 1  \, .
\]

Define $C' := e^{\frac{T_1}{2}}$ , let $U$ be a $C'$-rectangular neighbourhood of $1 \in G$ and let $\delta > 0$ such that $\overline{U_\delta (1)} \subset U$ . Then by the left invariance of the splitting and the metric on $G$ we see that $g U$ is a
$C'$-rectangular neighbourhood of $g$ and $\overline{U_\delta (g)} = g \overline{U_\delta (1)} \subset g U$ for all $g \in G$ . Define

\[
\eps_1 := \min\rund{\delta \frac{1 - e^{- \frac{T_1}{2}} }{e^{T_1} + 1}, \frac{T_1}{C_1}} > 0 \,  .
\]

Now assume $\gamma \in \Gamma$ and $T \geq T_1$ such that

\[
\eps := d\rund{\gamma x , x a_{T \bv}} \leq \eps_1 \, .
\]

Then $\varphi: G \rightarrow G \, , \, g \mapsto \gamma^{- 1} g a_T$ is a partially hyperbolic diffeomorphism with constant $e^{T_1} > 1$ and the corresponding splitting $T G = T^0 \oplus T^+ \oplus T^-$~. Then since

\[
\eps \leq \delta \frac{1 - e^{- \frac{T_1}{2}} }{e^{T_1} + 1} = \delta \frac{1 - C' e^{- T_1}}{C'^2 + 1}
\]

the partial {\sc Anosov} closing lemma, theorem \ref{closing} , tells us that there exist $y, z \in G$ such that

\begin{itemize}

\item[(i)] $x$ and $y$ belong to the same $T^-$-leaf and

\[
d^-\rund{x, y} \leq \eps \frac{C'}{1 - \frac{C'}{C}}  \, ,
\]

\item[(iii)] $y$ and $z$ belong to the same $T^+$-leaf and

\[
d^+\rund{y a_{T \bv} , z a_{T \bv}} \leq \eps \frac{C'^3}{1 - \frac{C'}{C}} \,  ,
\]

\item[(iv)] $\gamma z$ and $z a_{T \bv}$ belong to the same $T^0$-leaf and

\[
d^0\rund{\gamma z, z a_{T \bv}} \leq \eps C'^4 \, .
\]

\end{itemize}

In (iii) and (iv) we already used that the metric and the flow are left invariant. So by (iv) and since the $T^0$-leaf containing $z a_T$ is $z A M$~, there exist $w \in M$ and $t_0 \in \rz$ such that $\gamma z = z a_{t_0} w$ . So

\[
d^0\rund{a_{t_0 - T} w , 1} \leq \eps C'^4 \,  ,
\]

and so, since $A M \simeq \rz \times M$ isometrically, we see that

\[
d\rund{\rund{t_0, w} , \rund{T, 1}} \leq \eps C'^4 = \eps e^{2 T_1} \leq \eps C_1 \,  .
\]

In particular $\abs{t_0 - T} \leq T_1$ , and so $t_0 > 0$ . \\

Now let $\tau \in \, [ \, 0, T \, ]$ . Then since $x$ and $y$ belong to the same $T^-$-leaf the same is true for $x a_\tau$ and $y a_\tau$ , and

\[
d^-\rund{x a_\tau , y a_\tau} \leq d^-\rund{x, y} e^{- \tau} \leq \eps \frac{C'}{1 - \frac{C'}{C}} e^{- \tau} \leq \eps C_1 e^{- \tau} \,  .
\]

Since $y$ and $z$ belong to the same $T^+$-leaf the same is true for $y a_\tau$ and $z a_\tau$~, and

\begin{eqnarray*}
d^+\rund{y a_\tau , z a_\tau} &\leq& d^+\rund{y a_T , z a_T} e^{- \rund{T - \tau}} \\
&\leq& \eps \frac{C'^3}{1 - \frac{C'}{C}} e^{- \rund{T - \tau}} \leq \eps C_1 e^{- \rund{T - \tau}} \,  .
\end{eqnarray*}

Combining these two inequalities we obtain

\[
d\rund{x a_\tau , z a_\tau} \leq \eps C_1 \rund{e^{- \tau} + e^{- \rund{T - \tau}} } \, . \, \Box
\]

(ii) Let $L \subset G$ be compact and let $c \geq 1$ be given such that $\norm{\Ad_g} , \norm{\Ad_g^{- 1}} \leq c$ and therefore

\[
\frac{1}{c} \, d(a g, b g) \leq d(a, b) \leq c \, d(a g, b g)
\]

for all $g \in L$ and $a, b \in G$ . Let $\eps_0 > 0$ be given by lemma \ref{Gamma} and define

\[
\eps_2 := \frac{\eps_0}{c} > 0  \, .
\]

Let $x \in L$ , $\gamma \in \Gamma$ and $T \in \, \eckig{ \, 0, T_0 \, }$ such that

\[
\eps := d\rund{\gamma x, x a_T} \leq \eps_2  \, .
\]

Then since $x \in L$ we get

\[
d\rund{\gamma , x a_T x^{- 1}} \leq c \eps \leq \eps_0
\]

and so $\gamma \in \Gamma \cap N_{L, T_0, \eps_0}$ . This implies $\gamma = 1$ and so $d\rund{1, a_T} = \eps$ and therefore $T \leq \eps$ . $\Box$

\section{The unbounded realization} \label{unbounded realization}

Let $\n \sqsubset \g'$ be the standard maximal nilpotent sub {\sc Lie} algebra, which is at the same time the direct sum of all root spaces of $\g'$ of
positive roots with respect to $\a$ . Let $N := \exp \n$ . Then we have an {\sc Iwasawa} decomposition

\[
G = N A K  \, ,
\]

$N$ is $2$-step nilpotent, and so $N' := [N, N]$ is at the same time the center of~$N$~. \\

Now we transform the whole problem to the unbounded realization via the partial {\sc Cayley} transformation

\[
R := \rund{\begin{array}{c|c|c}
\frac{1}{\sqrt{2}} & 0 & \frac{1}{\sqrt{2}} \\ \hline
0 & 1 & 0 \\ \hline
- \frac{1}{\sqrt{2}} & 0 & \frac{1}{\sqrt{2}}
\end{array}}
\begin{array}{l}
\leftarrow 1 \\
\rbrace n - 1 \\
\leftarrow n + 1
\end{array} \in G'^\cz = SL(n + 1, \cz)
\]

mapping $B$ biholomorphically onto the unbounded domain

\[
H := \schweif{\left.\bw = \rund{\begin{array}{c}
w_1 \\ \hline
\bw_2
\end{array}} \begin{array}{l}
\leftarrow 1 \\
\rbrace n - 1
\end{array} \in \cz^n \, \right| \, \Re w_1 > \frac{1}{2} \bw_2^* \bw_2}  \, .
\]

We see that

\[
R G' R^{- 1} \sqsubset G'^\cz = SL(n + 1, \cz) \hookrightarrow GL(n + 1, \cz) \times GL(r, \cz)
\]

acts holomorphically and transitively on $H$ via fractional linear transformations, and explicit calculations show that

\[
a'_t := R a_t R^{- 1} = \rund{\begin{array}{c|c|c}
e^t & 0 & 0 \\ \hline
0 & 1 & 0 \\ \hline
0 & 0 & e^{- t}
\end{array}}
\begin{array}{l}
\leftarrow 1 \\
\rbrace n - 1 \\
\leftarrow n + 1
\end{array}
\]

for all $t \in \rz$ , and $R N R^{- 1}$ is the image of

\[
\rz \times \cz^{n - 1} \rightarrow R G' R^{- 1} \, , \, \rund{\lambda, \bu} \mapsto n'_{\lambda, \bu} := \rund{\begin{array}{c|c|c}
1 & \bu^* & i \lambda + \frac{1}{2} \bu^* \bu \\ \hline
0 & 1 & \bu \\ \hline
0 & 0 & 1
\end{array}}  \, ,
\]

which is a $\C^\infty$-diffeomorphism onto its image, with the multiplication rule

\[
n'_{\lambda, \bu} n'_{\mu, \bv} = n'_{\lambda + \mu + \Im\rund{\bu^* \bv} , \bu + \bv}
\]

for all $\lambda, \mu \in \rz$ and $\bu, \bv \in \cz^{n - 1}$ , so $N$ is exactly the {\sc Heisenberg} group $H_n$ acting on $H$ as pseudo translations

\[
\bw \mapsto \rund{\begin{array}{c}
w_1 + \bu^* \bw_2 + i \lambda + \frac{1}{2} \bu^* \bu \\ \hline
\bw_2 + \bu
\end{array}}  \, .
\]

Define $j\rund{R, \bz} =  \frac{\sqrt{2}}{1 - z_1} \in \O(B)$ , \\
$j\rund{R^{- 1}, \bw} := j\rund{R, R^{- 1} \bw}^{- 1} = \frac{\sqrt{2}}{1 + w_1} \in \O(H)$ , and for all

\[
g \in R G R^{- 1} = \rund{\begin{array}{c|c}
\begin{array}{c|c}
A & \bb \\ \hline
\bc & d
\end{array} & 0 \\ \hline
0 & E
\end{array}} \in R G R^{- 1}
\]

define

\[
j\rund{g, \bw} = j\rund{R, R^{- 1} g \bw} j\rund{R^{- 1} g R, R^{- 1} \bw} j\rund{R^{- 1}, \bw} = \frac{1}{\bc \bw + d} \, .
\]

Let $\H := H^{| r}$ with even coordinate functions $w_1, \dots, w_n$ and odd coordinate functions $\vartheta_1, \dots, \vartheta_r$ .  $R$ commutes with all $g \in Z_G\rund{G'}$ , and we have a right-representation of the group $R G R^{- 1}$ on $\D(\H)$ 
given by

\[
|_g : \D(\H) \rightarrow \D(\H) \, , \, f \mapsto f\rund{g \rund{\begin{array}{c}
\diamondsuit \\ \hline
\bthet
\end{array}} } j\rund{g, \diamondsuit}^k
\]

for all $g \in R G R^{- 1}$ . If we define

\[
|_R : \D(\H) \rightarrow \D(\B) \, , \, f \mapsto f\rund{R \rund{\begin{array}{c}
\diamondsuit \\ \hline
\bzeta
\end{array}} } j\rund{R, \diamondsuit}^k
\]

and

\[
|_{R^{- 1}} : \D(\B) \rightarrow \D(\H) \, , \, f \mapsto f\rund{R^{- 1} \rund{\begin{array}{c}
\diamondsuit \\ \hline
\bthet
\end{array}} } j\rund{R^{- 1}, \diamondsuit}^k  \, ,
\]

then we see that we get a commuting diagram

\[
\begin{array}{ccc}
\phantom{1234} \D(\H) & \mathop{\longrightarrow}\limits^{\phantom{1} |_{R g R^{- 1}} } & \D(\H) \\
|_R \phantom{1} \downarrow & \circlearrowleft & \phantom{12} \downarrow \phantom{1} |_R \\
\phantom{1234} \D(\B) & \mathop{\longrightarrow}\limits_{\phantom{1} |_g} & \D(\B)
\end{array}   \, .
\]

Now define the sesqui polynomial $\Delta'$ on $H \times H$ , holomorphic in the first and antiholomorphic in the second variable, as

\[
\Delta'\rund{\bz, \bw} := \Delta\rund{R^{- 1} \bz, R^{- 1} \bw} j\rund{R^{- 1}, \bz}^{- 1} \overline{j\rund{R^{- 1}, \bw}}^{- 1} = z_1 + \overline{w_1} - \bw_2^* \, \bz_2
\]

for all $\bz, \bw \in H$ . Clearly $\abs{\det\rund{\bz \mapsto R \bz}'} = \abs{j\rund{R, \bz}}^{n + 1}$ for all $\bz \in B$ . So

\[
\abs{\det\rund{\bw \mapsto g \bw}'} = \abs{j\rund{g, \bw}}^{n + 1} \, ,
\]

\[
\abs{j\rund{g, \be_1}} = \Delta'\rund{g \be_1, g \be_1}^{\frac{1}{2}}
\]

for all $g \in R G R^{- 1}$ and $\Delta'\rund{\bw, \bw}^{- (n + 1)} d V_\Leb$ is the $R G R^{- 1}$ -invariant volume element on $H$ . If
$f = \sum_{I \in \wp(r)} f_I \bzeta^I \in \O(\B)$ , all $f_I \in \O(B)^\cz$ , $I \in \wp(r)$~, then

\[
f|_{R^{- 1}} = \sum_{I \in \wp(r)} f_I \rund{R^{- 1} \diamondsuit} j\rund{R^{- 1}, \diamondsuit}^{k + \abs{I}} \bthet^I \in \O(\H)  \, ,
\]

and if $f = \sum_{I \in \wp(r)} f_I \bthet^I \in \O(\H)$ , all $f_I \in \C^\infty(H)^\cz$ , $I \in \wp(r)$ , and $g \in R G R^{- 1}$ , then

\[
f|_g = \sum_{I \in \wp(r)} f_I \rund{g \diamondsuit} j\rund{g, \diamondsuit}^{k + \abs{I}} \rund{E_g \bthet}^I \in \O(\H)  \, .
\]

Let $\partial H = \schweif{\bw \in \cz^n \, \left| \, \Re w_1 = \frac{1}{2} \bw_2^* \bw\right.}$ be the boundary of $H$ in $\cz^n$ . Then $\Delta'$ and $\partial H$ are $R N R^{- 1}$ -invariant, and
$R N R^{- 1}$ acts transitively on $\partial H$ and on each

\[
\schweif{\bw \in H \, \left| \, \Delta'\rund{\bw, \bw} = e^{2 t}\right.} = R N a_t \b0  \, ,
\]

$t \in \rz$ .

\pagebreak

\begin{figure}[H]
\begin{center}
\includegraphics[width=0.9\textwidth]{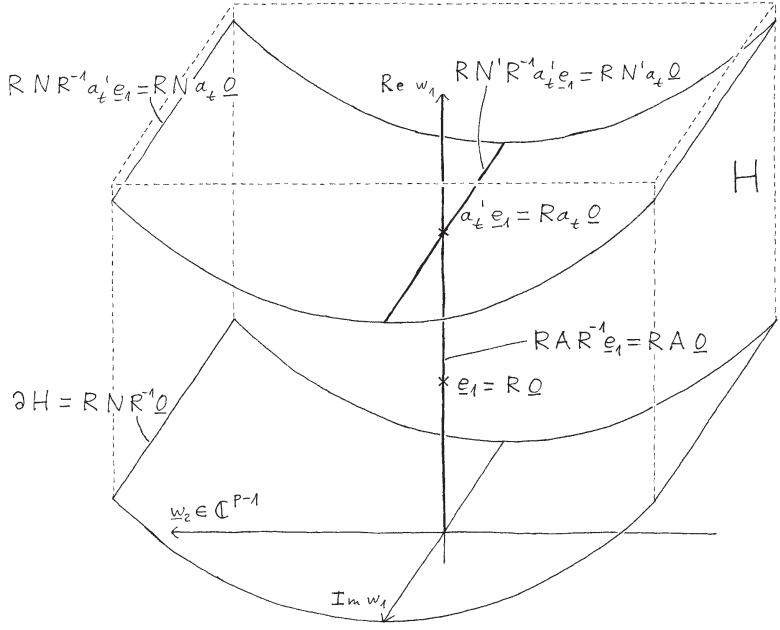}
\caption{the geometry of $H$ .}
\end{center}
\end{figure}

All geodesics in $H$ can be written in the form

\[
\rz \rightarrow H \, , \, t \mapsto \bw_t := R g a_t \b0 = R g R^{- 1} a'_t \be_1
\]

with some $g \in G$ , and conversely all these curves are geodesics in $H$ . We have to distinguish two cases: Either the goedesic connects $\infty$ with a point in $\partial H$ , or it
connects two points in $\partial H$ . In the second case we have

\[
\lim_{t \rightarrow \pm \infty} \Delta'\rund{\bw_t, \bw_t} = 0  \, ,
\]

so we may assume without loss of generality that $\Delta'\rund{\bw_t, \bw_t}$ is maximal for $t = 0$ , otherwise we have to reparametrize the geodesic using $g a_T$ , $T \in \rz$ appropiately
chosen, instead of $g$ .

\begin{lemma} \label{geodesic}
\item[(i)] Let

\[
\rz \rightarrow H \, , \, t \mapsto \bw_t := R g a_t \b0 = R g R^{- 1} a'_t \be_1
\]

be a geodesic in $H$ such that $\lim_{t \rightarrow \infty} \bw_t = \infty$ and $\lim_{t \rightarrow - \infty} \bw_t \in \partial H$ with respect to the euclidian metric on $\cz^p$ . Then for all $t \in \rz$

\[
\Delta'\rund{\bw_t, \bw_t} = e^{2 t} \Delta'\rund{\bw_0, \bw_0}  \, ,
\]

and if instead $\lim_{t \rightarrow - \infty} \bw_t = \infty$ and $\lim_{t \rightarrow \infty} \bw_t \in \partial H$ then

\[
\Delta'\rund{\bw_t, \bw_t} = e^{- 2 t} \Delta'\rund{\bw_0, \bw_0}  \, .
\]

\item[(ii)] Let

\[
\rz \rightarrow H \, , \, t \mapsto \bw_t := R g a_t \b0 = R g R^{- 1} a'_t \be_1
\]

be a geodesic in $H$ connecting two points in $\partial H$ such that $\Delta'\rund{\bw_t, \bw_t}$ is maximal for $t = 0$ . Then

\[
\rz \rightarrow \rz_{> 0} \, , \, t \mapsto \Delta'\rund{\bw_t, \bw_t}
\]

is strictly increasing on $\rz_{\leq 0}$ and strictly decreasing on $\rz_{\geq 0}$ , and for all $t \in \rz$

\[
\Delta'\rund{\bw_{- t}, \bw_{- t}} = \Delta'\rund{\bw_t, \bw_t}
\]

and

\[
e^{- 2 \abs{t}} \Delta'\rund{\bw_0, \bw_0} \leq \Delta'\rund{\bw_t, \bw_t} \leq 4 e^{- 2 \abs{t}} \Delta'\rund{\bw_0, \bw_0}  \, .
\]

\end{lemma}

{\it Proof:} (i) Since $R N R^{- 1}$ acts transitively on $\partial H$ and $\Delta'$ is $R N R^{- 1}$ -invariant we can assume without loss of generality that the geodesic connects $\b0$ and $\infty$ .
But in $H$ a geodesic is uniquely determined up to reparametrization by its endpoints. So we see that in the first case

\[
w_t = a'_t x \be_1 = e^{2 t} x \be_1
\]

and in the second case

\[
w_t = a'_{- t} x \be_1 = e^{- 2 t} x \be_1
\]

both with an appropriately chosen $x > 0$ . $\Box$  \\

(ii) Let $u, y \in \rz$ and $\bs \in \cz^{p - 1}$ such that $y^2 + \bs^* \bs = 1$ . Then

\[
\rz \rightarrow H \, , \, t \mapsto \bw^{\rund{u, y, \bs}}_t := \frac{e^u}{1 + y^2 \tanh^2 t} \rund{\begin{array}{c}
e^u \rund{1 - y^2 \tanh^2 t + 2 i y \tanh t} \\ \hline
\sqrt{2} \, \tanh t \rund{1 + i y \tanh t} \bs
\end{array}}
\]

is a geodesic through $e^{2 u} \be_1$ in $H$ since it is the image of the standard geodesic

\[
\rz \rightarrow B \, , \, t \mapsto a_t \b0 = \rund{\begin{array}{c}
\tanh t \\ \hline
\b0
\end{array}}
\]

in $B$ under the transformation

\[
\underbrace{\phantom{\begin{array}{c}
\\
\\
\\
\end{array}} a'_u \phantom{\begin{array}{c}
\\
\\
\\
\end{array}} }_{\in R A R^{- 1} \sqsubset R G' R^{- 1}} R \, \underbrace{\rund{\begin{array}{cc|c}
i y & - \bs^* & 0 \\
\bs & -i y & 0 \\ \hline
0 & 0 & 1
\end{array}} }_{\in K' \sqsubset G'}  \, .
\]

So we see that $\left.\partial_t \bw_t^{\rund{u, y, \bs}} \right|_{t = 0} = \rund{\begin{array}{c}
2 i e^{2 u} y \\ \hline
\sqrt{2} \, e^u \bs
\end{array}} \in T_{e^{2 u} \be_1} H$ is a unit vector with respect to the $R G R^{- 1}$ -invariant metric on $H$ coming from $B$ via $R$ . Now since $R N R^{- 1}$ acts transitively on each

\[
\schweif{\bw \in H \, \left| \, \Delta'\rund{\bw, \bw} = e^{2 t}\right.} = R N a_t \b0  \, ,
\]

$t \in \rz$ , and $\Delta'$ is invariant under $R N R^{- 1}$ we may assume without loss of generality that $\bw_0 = e^{2 u} \be_1$ with an appropriate $u \in \rz$ . Since $\Delta'\rund{\bw_t, \bw_t}$
is maximal for $t = 0$ we know that $\left.\partial_t \bw_t \right|_{t = 0}$ is a unit vector in\\
$i \rz \oplus \cz^{p - 1} \sqsubset T_{\be_1} H$ , and therefore there exist $y \in \rz$ and $\bs \in \cz^{p - 1}$ such that $y^2 + \bs^* \bs = 1$ and
$\left.\partial_t \bw_t \right|_{t = 0} = \rund{\begin{array}{c}
2 i e^{2 u} y \\ \hline
\sqrt{2} \, e^u \bs
\end{array}}$ . Since the geodesic is uniquely determined by $\bw_0$ and $\left.\partial_t \bw_t \right|_{t = 0}$ we see that $\bw_t = \bw^{\rund{u, y, \bs}}_t$ for all $t \in \rz$ , and so a straight
forward calculation shows that

\begin{eqnarray*}
\Delta'\rund{\bw_t, \bw_t} &=& 2 e^{2 u} \frac{1 - \tanh^2 t}{1 + y^2 \tanh^2 t} \\
&=& \frac{8 e^{ 2 u}}{\rund{1 + y^2}\rund{e^{2 t} + e^{- 2 t}} + 2 \bs^* \bs} \, .
\end{eqnarray*}

The rest is an easy exercise using $y^2 + \bs^* \bs = 1$ . $\Box$ \\

For all $t \in \rz$ define $A_{> t} := \schweif{\left.a_\tau \, \right| \, \tau > t} \subset A$ .

\begin{theorem}[a 'fundamental domain' for $\Gamma \backslash G$ ] \label{fund} There exist $\eta \subset N$ open and relatively
compact , $t_0 \in \rz$ and $\Xi \subset G'$ finite such that if we define

\[
\Omega := \bigcup_{g \in \Xi} g \eta A_{> t_0} K
\]

then

\item[(i)] $g^{- 1} \Gamma g \cap N Z_G\rund{G'} \sqsubset N Z_G\rund{G'}$ and $g^{- 1} \Gamma g \cap N' Z_G\rund{G'} \sqsubset N' Z_G\rund{G'}$ are lattices, and

\[
N Z_G\rund{G'} = \rund{g^{- 1} \Gamma g \cap N Z_G\rund{G'}} \eta Z_G\rund{G'}
\]

for all $g \in \Xi$ ,
\item[(ii)] $G = \Gamma \Omega$ ,
\item[(iii)] the set $\schweif{\gamma \in \Gamma \, | \, \gamma \Omega \cap \Omega \not= \emptyset}$ is finite.
\end{theorem}

{\it Proof:} direct consequence of theorem 0.6 (i) - (iii) , theorem 0.7 , lemma 3.16 and lemma 3.18 of \cite{GarlRagh} . For a detailed derivation see \cite{Kne} or section 3.2 of \cite{KneBuch} . $\Box$ \\

Now clearly the set of cusps of $\Gamma \backslash B$ in $\Gamma \backslash \partial B$ is contained in the set

\[
\schweif{\left.\lim_{t \rightarrow + \infty} \Gamma g a_t \b0 \, \right| \, g \in \Xi}  \, ,
\]

and is therefore finite as expected, where the limits are taken with respect to the Euclidian metric on $B$ .

\begin{cor} \label{L^s} Let $t_0 \in \rz$ , $\eta \subset N$ and $\Xi \subset G$ be given by theorem
\ref{fund}~. Let $h \in \C\rund{\Gamma \backslash G}^\cz$ and $s \in \, ] \, 0, \infty \, ] \,$ . Then $h \in L^s\rund{\Gamma \backslash G}$ if and only if
$h\rund{g \diamondsuit} \in L^s\rund{\eta A_{> t_0} K}$ for all $g \in \Xi$ .
\end{cor}

Let $f \in sM_k(\Gamma)$ and $g \in \Xi$ . Then we can decompose \\
$\left.f|_g \right|_{R^{- 1}} = \sum_{I \in \wp(r)} q_I \bthet^I \in \O(\H)$ , all $q_I \in \O(H)$ , $I \in \wp(r)$~, and by theorem \ref{fund} (i) we know that
$g^{- 1} \Gamma g \cap N' Z_G\rund{G'} \not\sqsubset Z_G\rund{G'}$ . So let \\
$n \in g^{- 1} \Gamma g \cap N' Z_G\rund{G'} \setminus Z_G\rund{G'}$ ,

\[
R n R^{- 1} = n'_{\lambda_0, \b0} \rund{\begin{array}{c|c}
\eps 1 & 0 \\ \hline
0 & E_0
\end{array}}  \, ,
\]

$\lambda_0 \in \rz \setminus \{0\}$ , $\eps \in U(1)$ , $E_0 \in U(r)$ , $\eps^{n + 1} = \det E$ . \\
$j\rund{R n R^{- 1}} := j\rund{R n R^{- 1}, \bw} = \eps^{- 1} \in U(1)$ is independent of $\bw \in H$ . So there exists $\chi \in \rz$ such that
$j\rund{R n R^{- 1}} = e^{2 \pi i \chi}$ . Without loss of generality we can assume that $E_0$ is diagonal, otherwise conjugate $n$ with an appropriate
element of $Z_G\rund{G'}$ . So there exists $D \in \rz^{r \times r}$ diagonal such that $E_0 = \exp\rund{2 \pi i D}$ .

\begin{theorem}[{\sc Fourier} expansion of $\left.f|_g \right|_{R^{- 1}}$ ] \label{Satake Fourier}

\item[(i)] There exist unique $c_{I, m} \in \O\rund{\cz^{n - 1}}$ , $I \in \wp(r)$ , \\
$m \in \frac{1}{\lambda_0} \rund{\zz - \tr_I D - \rund{k + \abs{I}} \chi}$ , such that

\[
q_I\rund{\bw} = \sum_{m \in \frac{1}{\lambda_0} \rund{\zz - \tr_I D - \rund{k + \abs{I}} \chi}} c_{I, m}\rund{\bw_2} e^{2 \pi m w_1}
\]

for all $\bw \in H$ and $I \in \wp(r)$ , and so

\[
\left.f|_g \right|_{R^{- 1}}\rund{\bw} = \sum_{I \in \wp(r)} \phantom{1} 
\sum_{m \in \frac{1}{\lambda_0} \rund{\zz - \tr_I D - \rund{k + \abs{I}} \chi}} c_{I, m}\rund{\bw_2} e^{2 \pi m w_1} \bthet^I
\]

for all $\bw = \rund{\begin{array}{c}
w_1 \\ \hline
\bw_2
\end{array}} \begin{array}{l}
\leftarrow 1 \\
\rbrace n - 1
\end{array} \in H$ , where the convergence is absolute and compact.

\item[(ii)] $c_{I, m} = 0$ for all $I \in \wp(r)$ and $m > 0$ (this is a super analogon for {\sc Koecher}'s principle, see for example in section 11.5 of \cite{Baily} ) , and if \\
$\tr_I D + \rund{k + \abs{I}} \chi \in \zz$ then $c_{I, 0}$ is a constant.

\item[(iii)] Let $I \in \wp(r)$ and $s \in \, [ \, 1, \infty \, ] \,$ . If $\tr_I D + \rund{k + \abs{I}} \chi \not\in \zz$ then

\[
q_I \Delta'\rund{\bw, \bw}^{\frac{k + \abs{I}}{2}} \in L^s\rund{R \eta A_{> t_0} \b0}
\]

with respect to the $R G R^{- 1}$ -invariant measure $\Delta'\rund{\bw, \bw}^{- (n + 1)} d V_\Leb$ on $H$~. If $\tr_I D + \rund{k + \abs{I}} \chi \in \zz$ and $k \geq 2 n - \abs{I}$ then

\[
q_I \Delta'\rund{\bw, \bw}^{\frac{k + \abs{I}}{2}} \in L^s\rund{R \eta A_{> t_0} \b0}
\]

with respect to the $R G R^{- 1}$ -invariant measure on $H$ if and only if $c_{I, 0} = 0$ .
\end{theorem}

A proof can be found in \cite{Kne} or \cite{KneBuch} section 3.2 .

\section{Proof of the main result}

We have a {\sc Lie} algebra embedding

\[
\rho: \sl(2, \cz) \hookrightarrow \g'^\cz = \sl(n + 1, \cz) \, , \, \rund{\begin{array}{cc}
a & b \\
c & - a
\end{array}} \mapsto \rund{\begin{array}{c|c}
\begin{array}{c|c}
a & 0 \\ \hline
0 & 0
\end{array} & \begin{array}{c}
b \\ \hline
0
\end{array} \\ \hline
\begin{array}{c|c}
c & 0
\end{array} & - a
\end{array}}   \, .
\]

Obviously the preimage of $\g'$ under $\rho$ is ${\frak su}(1, 1)$ , the preimage of $\k'$ under $\rho$ is ${\frak s}\rund{{\frak u}(1) \oplus {\frak u}(1)} \simeq {\frak u}(1)$ and $\rho$ lifts to a {\sc Lie} group homomorphism

\[
\widetilde \rho : SL(2, \cz) \rightarrow G'^\cz = SL(n + 1, \cz) \, , \, \rund{\begin{array}{cc}
a & b \\
c & d
\end{array}} \mapsto \rund{\begin{array}{c|c} \begin{array}{c|c}
a & 0 \\ \hline
0 & 0
\end{array} & \begin{array}{c}
b \\ \hline
0
\end{array} \\ \hline
\begin{array}{c|c}
c & 0
\end{array} & d
\end{array}}
\]

such that $\widetilde \rho \rund{SU(1,1)} \sqsubset G'$ .

Let us now identify the elements of $\g$ with the corresponding left invariant differential operators, they are defined on a dense subset of $L^2\rund{\Gamma \backslash G}$ , and define

\begin{eqnarray*}
&& \D := \rho\rund{\begin{array}{cc}
0&1 \\
1&0
\end{array}} \in \a \phantom{1} , \phantom{1} \D' := \rho\rund{\begin{array}{cc}
0&i \\
- i&0
\end{array}} \in \g'    \phantom{1} \text{and} \\
&& \phi := \rho\rund{\begin{array}{cc}
i&0 \\
0&- i
\end{array}} \in \k'  \, . \\
\end{eqnarray*}

The $\rz$-linear span of $\D$ , $\D'$ and $\phi$ is the $3$-dimensional sub {\sc Lie} algebra $\rho\rund{{\frak su} (1,1)}$ of $\g' \sqsubset \g$ , and $\D$ generates the flow $\varphi_t$ . $\phi$ generates a subgroup of $K'$ , being the image of the {\sc Lie}
group embedding

\[
\rz / 2 \pi \zz \hookrightarrow K \, , \, t \mapsto \exp\rund{t \phi} = \widetilde \rho \rund{\begin{array}{cc}
e^{i t} & 0 \\
0 & e^{- i t}
\end{array}} \, .
\]

Now define

\[
\D^+ := \frac{1}{2} \rund{\D - i \D'} \, , \, \D^- := \frac{1}{2} \rund{\D + i \D'} \text{  and  } \Psi := - i \phi
\]

as left invariant differential operators on $G$ . Then we get the commutation relations

\[
\eckig{\Psi, \D^+} = 2 \D^+ , \eckig{\Psi, \D^-} = - 2 \D^- \text{  and  }  \eckig{\D^+ , \D^-} = \Psi  \, ,
\]

and since $G$ is unimodular

\[
\rund{\D^+}^* = - \D^- , \rund{\D^-}^* = - \D^+ \text{ and } \Psi^* = \Psi  \, .
\]

So by standard {\sc Fourier} analysis

\[
L^2 \rund{\Gamma \backslash G} = \mathop{\widehat \bigoplus}\limits_{\nu \in \zz} H_\nu
\]

as an orthogonal sum, where

\[
H_\nu := \schweif{\left. F \in L^2\rund{\Gamma \backslash G} \cap \text{ domain } \Psi \, \right| \, \Psi F = \nu F}
\]

for all $\nu \in \zz$ . By a simple calculation we obtain

\[
\D^+ \rund{H_\nu \cap \text{ domain } \D^+} \subset H_{\nu + 2} \text{  and  } \D^- \rund{H_\nu \cap \text{ domain } \D^-} \subset H_{\nu - 2}
\]

for all $\nu \in \zz$ .

\begin{lemma} \label{holomorphic} $\D^- \widetilde h = 0$ for all $h \in \O(\B)$ . \end{lemma}

{\it Proof:} Let $g \in G$ . Then again $h |_g \in \O(\B)$ , and $\widetilde h \rund{g \diamondsuit} = \widetilde{h |_g}$ . So

\[
\D^- \widetilde h (g) = \D^- \rund{\widetilde h \rund{g \diamondsuit}} (1) = \overline{\partial}_1 h |_g (\b0) = 0 \, . \, \Box
\]

\begin{lemma}
Let $f \in sS_k^{(\rho)}(\Gamma)$ . Then $\widetilde f$ is uniformly {\sc Lipschitz} continuous.
\end{lemma}

{\it Proof:} Since on $G$ we use a left invariant metric it suffices to show that there exists a constant $c \geq 0$ such that for all $g \in G$ and $\xi \in \g$ with $\norm{\xi}_2 \leq 1$

\[
\abs{\xi \widetilde f(g)} \leq c \, .
\]

Then $c$ is a {\sc Lipschitz} constant for $\widetilde f$ . So choose an orthonormal basis $\rund{\xi_1, \dots, \xi_N}$ of $\g$ and a compact neighbourhood $L$ of $\b0$ in $B$ . Then by
{\sc Cauchy}'s integral formula there exist $C' , C'' \geq 0$ such that for all\\
$h \in \O(\B) \cap L_k^{\infty}(\B)$ and $n \in \{1, \dots, N\}$

\[
\abs{\rund{\xi_n \widetilde h}(1)} \leq C' \int_L \abs{h} \leq C' \vol L \norm{h}_{\infty, L} \leq C'' \vol L \norm{\widetilde h}_\infty \,  ,
\]

and since $\g \rightarrow \cz \, , \, \xi \mapsto \rund{\xi \widetilde h}(1)$ is linear we obtain

\[
\abs{\rund{\xi \widetilde h}(1)} \leq N C'' \vol L \norm{\widetilde h}_\infty
\]

for general $\xi \in \g$ with $\norm{\xi}_2 \leq 1$ . Now let $g \in G$ . Then again $f |_g \in \O(\B)$~, $\widetilde f \rund{g \diamondsuit} = \widetilde{f |_g}$ , and by {\sc Satake}'s theorem, theorem \ref{Satake} , $f$ and so \\
$f|_g \in L_k^\infty (\B)$ . So

\[
\abs{\xi \widetilde f(g)} = \abs{\rund{\xi \widetilde f \rund{g \diamondsuit}}(1)} \leq N C'' \vol L \norm{\widetilde f \rund{g \diamondsuit}}_\infty \leq N C'' \vol L \norm{\widetilde f}_\infty  \, ,
\]

and we can define $c := N C'' \vol L \norm{\widetilde f}_\infty$ . $\Box$ \\

Now let $f \in sS_k^{(\rho)}(\Gamma)$ such that $\rund{\varphi_{\gamma_0 , I, m}, f}_\Gamma = 0$ for all $\varphi_{\gamma_0 , I, m}$ , $\gamma_0 \in \Gamma$ primitive loxodromic, $I \in \wp(r)$ , $\abs{I} = \rho$ ,
$m \in \, ] \, - C, C \, [$ . We will show that $f = 0$ in several steps.

\begin{lemma} \label{Lipschitz}
There exists $F \in \C \rund{\Gamma \backslash G}^\cz \otimes \bigwedge\rund{\cz^r}$ uniformly {\sc Lipschitz} continuous on compact sets and differentiable along the flow $\varphi_t$ such that

\[
f = \left. \partial_\tau F\rund{\diamondsuit a_\tau} \right|_{\tau = 0} = \D F \, .
\]
\end{lemma}

{\it Proof:} Here we use that the right translation with $A$ on $\Gamma \backslash G' \widetilde Z$ is topologically transitive. So let $g_0 \in G'$ such that $\overline{\Gamma g_0 A} = G' \widetilde Z$ and define 
$s \in \C^\infty\rund{\rz}^\cz \otimes \bigwedge\rund{\cz^r}$ by

\[
s(t) := \int_0^t \widetilde f\rund{g_0 a_\tau} d \tau
\]

for all $t \in \rz$ . \\

{\it Step I} {\bf Show that for all $L \subset G' \widetilde Z$ compact there exist constants $C_3 \geq 0$ and $\eps_3 > 0$ such that for all $t \in \rz$ , $T \geq 0$ and $\gamma \in \Gamma$ if
$g_0 a_t \in L$ and

\[
\eps := d\rund{\gamma g_0 a_t , g_0 a_{t + T}} \leq \eps_3
\]

then $\abs{s(t) - s(t + T)} \leq C_3 \eps$ . } \\

Let $L \subset G' \widetilde Z$ be compact, $T_0 > 0$ be given by lemma \ref{Gamma} and $C_1 \geq 1$ and $\eps_1$ be given by theorem \ref{Anosov} (i) with
$T_1 := T_0$ . Define $C_3 := \max\rund{C_1 \rund{C_2 + 2 c}, \norm{\widetilde f}_\infty} \geq 0$ , where $C_2 \geq 0$ is the {\sc Lipschitz} constant from theorem \ref{Fourier} (ii)
and $c \geq 0$ is the {\sc Lipschitz} constant of $\widetilde f$~. Define $\eps_3 := \min\rund{\eps_1, \eps_2, \frac{T_0}{2 C_1}} > 0$ , where $\eps_2 > 0$ is given by theorem
\ref{Anosov}~(ii)~. \\

Let $t \in \rz$ , $T \geq 0$ and $\gamma \in \Gamma$ such that $g_0 a_t \in L$ and $\eps := d\rund{\gamma g_0 a_t , g_0 a_{t + T}} \leq \eps_3$ . \\

First assume $T \geq T_0$ . Then by theorem \ref{Anosov} (i) since $\eps \leq \eps_1$ there exist \\
$g \in G$ , $w_0 \in M$ and $t_0 > 0$ such that $\gamma g = g a_{t_0} w_0$ , \\
$d\rund{\rund{t_0, w_0} , \rund{T, 1}} \leq C_1  \eps$~, and for all $\tau \in \, [ \, 0, T \, ]$

\[
d\rund{g_0 a_{t + \tau} , g a_\tau} \leq C_1 \eps \rund{e^{- \tau} + e^{- \rund{T - \tau}} } \,  .
\]

We get

\[
s(t + T) - s(t) = \underbrace{\int_0^T \widetilde f \rund{g a_\tau} d \tau}_{I_1 :=}
+ \underbrace{\int_0^T \rund{\widetilde f \rund{g_0 a_{t + \tau}} - \widetilde f \rund{g a_\tau}} d \tau}_{I_2 :=}
\]

and

\begin{eqnarray*}
\abs{I_2} &\leq& \int_0^T \abs{\widetilde f \rund{g_0 a_{t + \tau} } - \widetilde f \rund{g a_\tau} } d \tau \\
&\leq& c \int_0^T d\rund{g_0 a_{t + \tau} , g a_\tau} d \tau \\
&\leq& c C_1 \eps \int_0^T \rund{e^{- \tau} + e^{- \rund{T - \tau}} } d \tau \\
&\leq& 2 c C_1 \eps \,  .
\end{eqnarray*}
 \\

Since $\gamma \in \Gamma$ is regular loxodromic there exists $\gamma_0 \in \Gamma$ primitive loxodromic and $\nu \in \nz \setminus \{0\}$ such that $\gamma = \gamma_0^\nu$ . $\gamma_0 \in g A W g^{- 1}$ since lemma \ref{determined} tells us that 
$g \in G$ is already determined by $\gamma$ up to right translation with elements of $A N_K(A)$ . Choose $w' \in N_K(M)$ , $t_0' > 0$ and $w_0' \in M$ such that $E_{w_0'}$ is diagonal and $\gamma = g w' a_{t_0'} w_0' \rund{g w'}^{- 1}$ , and let 
$g' := g w'$ . We define $h \in \C^\infty\rund{\rz \times M}^\cz \otimes \bigwedge\rund{\cz^r}$ as

\[
h(\tau , w) := \widetilde f\rund{g' a_\tau w} = \widetilde f\rund{g a_\tau w' w}
\]

for all $\tau \in \rz$ and $w \in M$ . Then

\[
I_1 = \int_0^T h\rund{\tau, w'^{- 1}} d \tau \, .
\]

We can apply theorem \ref{Fourier} (i) and, since $f$ is perpendicular to all $\varphi_{\gamma_0 , I, m}$ , $I \in \wp(r)$ , $m \in \, ] \, - C, C \, [$ , also \ref{Fourier} (ii) with $g' := g w'$ instead of $g$ , and so

\begin{eqnarray*}
\abs{I_1} &=& \abs{H\rund{T, w'^{- 1}} - H\rund{0, w'^{- 1}} } \\
&=& \abs{H\rund{T, w'^{- 1}} - H\rund{t_0, w'^{- 1} w_0} } \\
&\leq& C_2 d\rund{\rund{T, 1} , \rund{t_0, w_0}} \\
&\leq& C_1 C_2 \eps  \, ,
\end{eqnarray*}

where we used that $H\rund{0, w'^{- 1}} = H\rund{t_0', w_0' w'^{- 1}}$ and that we have chosen the left invariant metric on $M$ , and the claim follows.

Now assume $T \leq T_0$ . Then by theorem \ref{Anosov} (ii) since $\eps \leq \eps_0$ we get $T \leq \eps$ and so

\[
\abs{s(t + T) - s(t)} = \abs{\int_0^T \widetilde f\rund{g_0 a_{t + \tau}} d \tau} \leq \eps \norm{\widetilde f}_\infty  \, .
\]

{\it Step II} {\bf Show that there exists a unique $F_1 \in \C\rund{\Gamma \backslash G' \widetilde Z}^\cz \otimes \bigwedge\rund{\cz^r}$ uniformly {\sc Lipschitz} continuous on compact sets such that for all $t \in \rz$

\[
s(t) = F_1\rund{g_0 a_t} \, .
\]  }

By step I for all $L \subset \Gamma \backslash G' \widetilde Z$ compact with $L^\circ \mathop{\subset}\limits_{\rm{dense}} L$ there exists a unique $F_L \in \C\rund{\Gamma \backslash G' \widetilde Z}^\cz$ uniformly {\sc Lipschitz} continuous such that for
all $t \in \rz$ if $\Gamma g_0 a_t \in L$ then $s(t) = F_L \rund{\Gamma g_0 a_t}$ . So we see that there exists a unique $F_1 \in \C\rund{\Gamma \backslash G' \widetilde Z}^\cz \otimes \bigwedge\rund{\cz^r}$ such that $\left.F_1 \right|_L = F_L$ for all
$L \subset \Gamma \backslash G' \widetilde Z$ compact with $L^\circ \mathop{\subset}\limits_{\text{dense}} L$~. \\

{\it Step III} {\bf Show that $F_1$ is differentiable along the flow and that for all $g \in G' \widetilde Z$

\[
\partial_\tau F_1\rund{g a_\tau} |_{\tau = 0} = \widetilde f (g)  \, .
\]  }

Let $g \in G' \widetilde Z$ . It suffices to show that for all $T \in \rz$

\[
\int_0^T \widetilde f\rund{g a_\tau} d \tau = F_1\rund{g a_T} - F_1(g) \, .
\]

If $g = g_0 a_t$ for some $t \in \rz$ then it is clear by construction. For general $g \in G' \widetilde Z$ since $\overline{\Gamma g_0 A} = G' \widetilde Z$ there exists $\rund{\gamma_n, t_n}_{n \in \nz} \in \rund{\Gamma \times \rz}^\nz$ such that

\[
\lim_{n \rightarrow \infty} \gamma_n g_0 a_{t_n} = g  \, ,
\]

and so

\[
\lim_{n \rightarrow \infty} \gamma_n g_0 a_{\tau + t_n} = g a_\tau
\]

compact in $\tau \in \rz$ , finally $\widetilde f$ is uniformly {\sc Lipschitz} continuous. Therefore we can interchange integration and taking limit $n \leadsto \infty$ :

\begin{eqnarray*}
\int_0^T \widetilde f\rund{g a_\tau} d \tau &=& \lim_{n \rightarrow \infty} \int_0^T \widetilde f\rund{\gamma_n g_0 a_{\tau + t_n}} d \tau \\
&=& \lim_{n \rightarrow \infty} \rund{F_1\rund{\gamma_n g_0 a_{T + t_n}} - F_1\rund{\gamma_n g_0 a_{t_n}} } \\
&=& F_1\rund{g a_T} - F_1(g) \, .
\end{eqnarray*}

{\it Step IV} {\bf Conclusion.} \\

Define $F \in \C(G)^\cz \otimes \bigwedge\rund{\cz^r}$ as

\[
F\rund{g w} := \int_{\widetilde Z} F_1\rund{g u^{- 1}, E_{u w} \bEta} j(u w)^{k + \rho} d u
\]

for all $g \in G' \widetilde Z$ and $w \in Z_G\rund{G'}$ , where we normalize the {\sc Haar} measure on the compact {\sc Lie} group $\widetilde Z$ such that $\vol \widetilde Z = 1$ . Then we see that $F$ is well defined and fulfills all the desired properties.
$\Box$ \\

\pagebreak

\begin{lemma} \label{differentiable}
\item[(i)] For all $L \subset G$ compact there exists $\eps_4 > 0$ such that for all $g, h \in L$ if $g$ and $h$ belong to the same $T^-$-leaf and $d^-(g, h) \leq \eps_4$ then

\[
\lim_{t \rightarrow \infty} \rund{F\rund{g a_t} - F\rund{h a_t}} = 0  \,  ,
\]

and if $g$ and $h$ belong to the same $T^+$-leaf and $d^+(g, h) \leq \eps_4$ then

\[
\lim_{t \rightarrow - \infty} \rund{F\rund{g a_t} - F\rund{h a_t}} = 0  \, .
\]

\item[(ii)] $F$ is continuously differentiable along $T^-$- and $T^+$-leafs, more precisely if $\rho: I \rightarrow G$ is a continuously differentiable curve in a $T^-$-leaf then

\[
\partial_t \rund{F \circ \rho} (t) = - \int_0^\infty \partial_t \widetilde f \rund{\rho(t) a_\tau} d \tau   \,  ,
\]

and if $\rho: I \rightarrow G$ is a continuously differentiable curve in a $T^+$-leaf then

\[
\partial_t \rund{F \circ \rho} (t) = \int_{- \infty}^0 \partial_t \widetilde f \rund{\rho(t) a_\tau} d \tau   \,  .
\]

\end{lemma}

{\it Proof:} (i) Let $L \subset G$ be compact, and let $L' \subset G$ be a compact neighbourhood of $L$ . Let $T_0 > 0$ be given by lemma \ref{Gamma} and $\eps_2 > 0$ by
theorem \ref{Anosov} (ii) both with respect to $L'$ . Define

\[
\eps_4 := \frac{1}{3} \min\rund{\eps_1, \eps_2, \frac{T_0}{2 C_1}} > 0 \, ,
\]

where $\eps_1 > 0$ and $C_1 \geq 1$ are given by theorem \ref{Anosov} (i) with $T_1 := T_0$~. Let $\delta_0 > 0$ such that $\overline{U_{\delta_0} (L)} \subset L'$ and let

\[
\delta \in \, \left] \, 0, \min\rund{\delta_0, \eps_4} \, \right[  \, .
\]

Let $g, h \in L$ be in the same $T^-$-leaf such that $\eps := d^-(g, h) \leq \eps_4$ . Since the splitting of $T G$ is left invariant and $T_1^{-} (G) \sqsubset \g'$ we see that there exist $g', h' \in G'$ and $u \in Z_G\rund{G'}$ such that $g = g' u$ and $h = h' u$ . 
Fix some $T' > 0$~. Again by assumption there exists $g_0 \in G'$ such that $\overline{\Gamma g_0 A} = G' \widetilde Z$ , and so $g, h \in \overline{\Gamma g_0 u A}$ . So there exist $\gamma_g , \gamma_h \in \Gamma$ and $t_g , t_h \in \rz$ such that

\[
d\rund{g a_t , \gamma_g g_0 u a_{t_g + t}} , d\rund{h a_t , \gamma_h g_0 u a_{t_h + t}} \leq \delta
\]

for all $t \in \, [ \, 0, T' \, ]$ , and so in particular $\gamma_g g_0 u a_{t_g} , \gamma_h g_0 u a_{t_h} \in L'$ . We will show that for all $t \in \, [ \, 0, T' \, ]$

\[
\abs{F\rund{\gamma_g g_0 u a_{t_g + t}} - F\rund{\gamma_h g_0 u a_{t_h + t}} } \leq C_3' \rund{\eps e^{- t} + 2 \delta}
\]

with the same constant $C_3' \geq 0$ as in step I of the proof of lemma \ref{Lipschitz} with respect to $L'$ .

\begin{quote}
Without loss of generality we may assume $T := t_h - t_g \geq 0$ . Define $\gamma := \gamma_g \gamma_h^{- 1} \in \Gamma$ . Then for all $t \in \, [ \, 0, T' \, ]$

\[
d\rund{\gamma \gamma_h g_0 u a_{t_g + t} , \gamma_h g_0 u a_{t_g + t + T}} \leq \eps e^{- t} + 2 \delta \, .
\]

First assume $T \geq T_0$ and fix $t \in \, [ \, 0, T' \, ]$ . Then by theorem \ref{Anosov} (i) since $\eps e^{- t} + 2 \delta \leq \eps + 2 \delta \leq \min\rund{\eps_1, \frac{T_0}{2 C_1}}$ there exist $z \in G$ , $t_0 \in \rz$ and $w \in M$ such that 
$\gamma z = z a_{t_0} w$ ,

\[
d\rund{\rund{t_0, w}, (T, 1)} \leq C_1 \rund{2 \delta + \eps e^{- t}} \, ,
\]

and for all $\tau \in \, [ \, 0, T \, ]$

\[
d\rund{\gamma_g g_0 u a_{t_g + t + \tau} , z a_\tau} \leq C_1 \rund{\eps e^{- t} + 2 \delta} \rund{e^{- \tau} + e^{- \rund{T - \tau}} }  \, .
\]

And so by the same calculations as in the proof of lemma \ref{Lipschitz} we obtain the estimate

\[
\abs{F\rund{\gamma_g g_0 u a_{t_g + t}} - F\rund{\gamma_h g_0 u a_{t_h + t}} } \leq C_3' \rund{\eps e^{- t} + 2 \delta}  \,  .
\]

Now assume $T \leq T_0$ . Then by theorem \ref{Anosov} (ii) since \\
$\gamma_g g_0m  a_{t_g} \in L'$ and $\eps + 2 \delta \leq \eps_2$ we obtain $\gamma = 1$ and so by the left invariance of the metric on $G$

\[
d\rund{1 , a_T} \leq \eps e^{- T'} + 2 \delta  \, ,
\]

therefore $T \leq \eps e^{- T'} + 2 \delta$ . So as in the proof of lemma \ref{Lipschitz}

\begin{eqnarray*}
\abs{F\rund{\gamma_g g_0 u a_{t_g + t}} - F\rund{\gamma_h g_0 u a_{t_h + t}} } &\leq& \norm{\widetilde f}_\infty \rund{\eps e^{- T'} + 2 \delta} \\
&\leq& C_3' \rund{\eps e^{- t} + 2 \delta}  \, .
\end{eqnarray*}

\end{quote}

Now let us take the limit $\delta \leadsto 0$ . Then $\gamma_g g_0 u a_{t_g} \leadsto g$ and $\gamma_h g_0 u a_{t_h} \leadsto h$ , so since $F$ is continuous

\[
\abs{F\rund{g a_t} - F\rund{h a_t}} \leq C_3' \eps e^{- t}
\]

for all $t \in [0, T']$ , and since $T' > 0$ has been arbitrary, we obtain this estimate for all $t \geq 0$ and so $\lim_{t \rightarrow \infty} F\rund{g a_t} - F\rund{h a_t} = 0$ . By similar
calculations we can prove $\lim_{t \rightarrow - \infty} F\rund{g a_t} - F\rund{h a_t} = 0$ if $g$ and $h$ belong to the same $T^+$-leaf and $d^+\rund{g, h} \leq \eps_4$ . $\Box$ \\

(ii) Let $\rho: I \rightarrow G$ be a continuously differentiable curve in a $T^-$-leaf, and let $t_0, t_1 \in I$ , $t_1 > t_0$ . It suffices to show that

\[
F\rund{\rho\rund{t_1}} - F\rund{\rho\rund{t_0}} = - \int_{t_0}^{t_1} \int_0^\infty \partial_t \widetilde f \rund{\rho(t) a_\tau} d \tau d t \, .
\]

Let $C' \geq 0$ such that $\norm{\partial_t \rho(t)} \leq C'$ for all $t \in \eckig{t_0, t_1}$ . Then since $\rho$ lies in a $T^-$-leaf we have $\norm{\partial_t \rund{\rho(t) a_\tau}} \leq C' e^{- \tau}$
and so

\[
\abs{\partial_t \widetilde f \rund{\rho(t) a_\tau}} \leq c C' e^{- \tau}
\]

for all $\tau \geq 0$ and $t \in \eckig{t_0, t_1}$ where $c \geq 0$ is the {\sc Lipschitz} constant of $\widetilde f$ . So the double integral on the right side is absolutely convergent and so we can
interchange the order of integration:

\begin{eqnarray*}
\int_{t_0}^{t_1} \int_0^\infty \partial_t \widetilde f \rund{\rho(t) a_\tau} d \tau d t &=& \int_0^\infty \int_{t_0}^{t_1} \partial_t \widetilde f \rund{\rho(t) a_\tau} d t d \tau \\
&=& \int_0^\infty \rund{\widetilde f \rund{\rho\rund{t_1} a_\tau} - \widetilde f \rund{\rho\rund{t_0} a_\tau}} d \tau \\
&=& \lim_{T \rightarrow \infty} \rund{F \rund{\rho\rund{t_1} a_T} - F \rund{\rho\rund{t_0} a_T}} \\
&& - F \rund{\rho\rund{t_1}} + F \rund{\rho\rund{t_0}} \,  .
\end{eqnarray*}

Now let $L \subset G$ be compact such that $\rho(\eckig{t_1, t_2}) \subset L$ and let $\eps_4 > 0$ as in (i) . Without loss of generality we may assume that
$d^-\rund{\rho\rund{t_0}, \rho\rund{t_1}} \leq \eps_4$~. Then

\[
\lim_{T \rightarrow \infty} \rund{F \rund{\rho\rund{t_1} a_T} - F \rund{\rho\rund{t_0} a_T}} = 0
\]

by (i) . By similar calculations one can also prove 

\[
\partial_t \rund{F \circ \rho} (t) = \int_{- \infty}^0 \partial_t \widetilde f \rund{\rho(t) a_\tau} d \tau
\]

in the case when $\rho: I \rightarrow G$ is a continuously differentiable curve in a $T^+$-leaf. $\Box$ \\

\begin{lemma} \label{L^2}
\item[(i)] $F \in L^2\rund{\Gamma \backslash G} \otimes \bigwedge\rund{\cz^r}$ ,
\item[(ii)] $\xi F \in L^2\rund{\Gamma \backslash G} \otimes \bigwedge\rund{\cz^r}$ for all $\xi \in \rz \D \oplus \g \cap \rund{T^+ \oplus T^-}$ .
\end{lemma}

{\it Proof:} (i) If $\Gamma \backslash G$ is compact then the assertion is trivial. So assume that $\Gamma \backslash G$ is not compact, then we use the unbounded realization $\H$ of $\B$ introduced in section \ref{unbounded realization} . Since
$\vol (\Gamma \backslash G) < \infty$ it suffices to prove that $F$ is bounded, and by corollary \ref{L^s} it is even enough to show that $F\rund{g \diamondsuit}$ is bounded on $N A_{> t_0} K$ for all $g \in \Xi$ , where $t_0 \in \rz$ and $\Xi \subset G'$ are 
given by theorem \ref{fund} . So let $g \in \Xi$ . \\

{\it Step I} {\bf Show that $F\rund{g \diamondsuit}$ is bounded on $N a_{t_0} K$ .} \\

Let also $\eta \subset N$ be given by theorem \ref{fund} . Then $F\rund{g \diamondsuit}$ is clearly bounded on the compact set $\overline \eta a_{t_0} K$ . On the other hand $F\rund{g \diamondsuit}$ is left- $g^{- 1} \Gamma g$ -invariant, so it is also
bounded on

\[
N a_{t_0} K = \rund{g \Gamma g^{- 1} \cap N Z_G\rund{G'}} \eta a_{t_0} K
\]

by theorem \ref{fund} (i) . \\

{\it Step II} {\bf Show that there exists $C' \geq 0$ such that for all $g' \in N A_{> t_0} K$

\[
\abs{\widetilde f\rund{g g'}} \leq \frac{C'}{\Delta'\rund{R g' \b0, R g' \b0}}  \, .
\]} \\

As in section \ref{unbounded realization} let $q_I \in \O(H)$ such that $\left.f|_g \right|_{R^{- 1}} = \sum_{I \in \wp(r)} q_I \bthet^I$ . Then since
$\widetilde f\rund{g \diamondsuit} \in L^2\rund{\eta A_{> t_0} K} \otimes \bigwedge\rund{\cz^r}$ by theorem \ref{Satake Fourier} we have {\sc Fourier} expansions

\begin{equation}
q_I\rund{\bw} = \sum_{m \in \frac{1}{\lambda_0} \rund{\zz - \tr_I D - \rund{k + \abs{I}} \chi} \cap \rz_{< 0}} c_{I, m}\rund{\bw_2} e^{2 \pi m w_1}  \label{eq Satake Fourier}
\end{equation}

for all $I \in \wp(I)$ and $\bw = \rund{\begin{array}{c}
w_1 \\ \hline
\bw_2
\end{array}} \begin{array}{l}
\leftarrow 1 \\
\rbrace n - 1
\end{array} \in H$ , where $c_{I, m} \in \O\rund{\cz^{n - 1}}$ , $I \in \wp(r)$ , $m \in \frac{1}{\lambda_0} \rund{\bz - \tr_I D - \rund{k + \abs{I}} \chi} \cap \rz_{< 0}$ . Define

\[
M_0 := \max \bigcup_{I \in \wp(r)} \frac{1}{\lambda_0} \rund{\zz - \tr_I D - \rund{k + \abs{I}} \chi} \cap \rz_{< 0} < 0 \, .
\]

$R \overline \eta a_{t_0} \b0 \subset H$ is compact, and so since the convergence of the {\sc Fourier} series (\ref{eq Satake Fourier}) is absolute and compact we can define

\begin{eqnarray*}
&& C'' := e^{- 2 \pi M_0 e^{2 t_0}} \times \\
&& \phantom{1} \times \max_{I \in \wp(r)} \sum_{m \in \frac{1}{\lambda_0} \rund{\zz - \tr_I D - \rund{k + \abs{I}} \chi} \cap \rz_{< 0}} \norm{c_{I, m}\rund{\bw_2} e^{2 \pi m w_1}}_{\infty, R \overline \eta a_{t_0} \b0} < \infty  \, .
\end{eqnarray*}

Then we have

\[
\abs{q_I\rund{\bw}} \leq C'' e^{\pi M_0 \Delta'\rund{\bw, \bw}}
\]

for all $I \in \wp(r)$ and $\bw \in R \eta A_{> t_0} \b0$ . Now let $g' = \rund{\begin{array}{c|c}
* & 0 \\ \hline
0 & E'
\end{array}} \in \eta A_{> 0} K$ , $E' \in U(r)$ . Then

\begin{eqnarray*}
\widetilde f\rund{g g'} &=& \left.\left. f|_g \right|_{R^{- 1}} \right|_{R g R^{- 1}} \rund{\be_1} \\
&=& \left.f|_g \right|_{R^{- 1}}\rund{R g' R^{- 1} \rund{\begin{array}{c}
\be_1 \\ \hline
\bEta
\end{array}} } j\rund{R g' R^{- 1}, \be_1}^k \\
&=& \left.f|_g \right|_{R^{- 1}}\rund{\begin{array}{c}
R g' \b0 \\ \hline
E \bEta j\rund{R g' R^{- 1}}
\end{array}} j\rund{R g' R^{- 1}, \be_1}^k \\
&=& \sum_{I \in \wp(r)} q_I\rund{R g' \b0} \rund{E \bEta}^I j\rund{R g' R^{- 1}, \be_1}^{k + \abs{I}}  \, .
\end{eqnarray*}

Therefore since $\abs{j\rund{R g' R^{- 1}, \be_1}} = \sqrt{\Delta'\rund{R g' \b0, R g' \b0}}$ we get

\begin{eqnarray*}
\abs{\widetilde f\rund{g g'}} &\leq& 2^r C'' e^{\pi M_0 \Delta'\rund{R g' \b0, R g' \b0}} \times \\
&& \phantom{1} \times \rund{\Delta'\rund{R g' \b0, R g' \b0}^{\frac{k}{2}} + \Delta'\rund{R g' \b0, R g' \b0}^{\frac{k + r}{2}} }  \, .
\end{eqnarray*}

So we see that there exists $C' > 0$ such that

\[
\abs{\widetilde f\rund{g g'}} \leq \frac{C'}{\Delta'\rund{R g' \b0, R g' \b0}}
\]

for all $g' \in \eta A_{> t_0} K$ , but on one hand $\widetilde f\rund{g \diamondsuit}$ is left- $g^{- 1} \Gamma g$ -invariant, and on the other hand $\Delta'$ is $R N Z_G\rund{G'} R^{- 1}$
-invariant. Therefore the estimate is correct even for all

\[
g' \in N A_{> t_0} K = \rund{g \Gamma g^{- 1} \cap N Z_G\rund{G'}} \eta A_{> t_0} K
\]

by theorem \ref{fund} (i) . \\

{\it Step III} {\bf Conclusion: Prove that

\[
\abs{F\rund{g \diamondsuit}} \leq \norm{F\rund{g \diamondsuit}}_{\infty, N a_{t_0} K} + 2 C' e^{- 2 t_0}
\]

on $N A_{> t_0} K$ .} \\

Let $g' \in G$ be arbitrary. We will show the estimate on $g' A \cap N A_{> t_0} K$ .

\[
\rz \rightarrow H \, , \, t \mapsto \bw_t := R g' a_t \b0
\]

is a geodesic in $H$ , and for all $t \in \rz$ we have $g' a_t \in N A_{> t_0} K$ if and only if $\Delta'\rund{\bw_t, \bw_t} > 2 e^{2 t_0}$ . Now we have to distinguish two cases. \\

In the first case the geodesic connects $\infty$ with a point in $\partial H$ . First assume that $\lim_{t \rightarrow \infty} \bw_t = \infty$ and $\lim_{t \rightarrow - \infty} \bw_t \in \partial H$ . Then
$\lim_{t \rightarrow \infty} \Delta'\rund{\bw_t, \bw_t} =~\infty$ and $\lim_{t \rightarrow - \infty} \Delta'\rund{\bw_t, \bw_t} = 0$ . So we may assume without loss of generality that
$\Delta'\rund{\bw_0, \bw_0} = 2 e^{2 t_0}$ , and therefore $g' = g' a_0 \in N a_{t_0} K$ and $g' a_t \in N A_{> t_0} K$ if and only if $t > 0$ . So let $t > 0$ . Then

\[
F\rund{g g' a_t} = F\rund{g g'} + \int_{0}^{t} \widetilde f\rund{g g' a_\tau} d \tau \, ,
\]

and so

\[
\abs{F\rund{g g' a_t}} \leq \norm{F\rund{g \diamondsuit}}_{\infty, N a_{t_0} K} + \int_{0}^{t} \abs{\widetilde f\rund{g g' a_\tau}} d \tau  \, .
\]

By step II and lemma \ref{geodesic} (i)

\begin{eqnarray*}
\int_{0}^{t} \abs{\widetilde f\rund{g g' a_\tau}} d \tau &\leq& C' \int_{0}^{t} \frac{d \tau}{\Delta'\rund{\bw_\tau, \bw_\tau}} \\
&=& \frac{C'}{\Delta'\rund{\bw_0, \bw_0}} \int_{0}^{t} e^{- 2 \tau} d \tau \\
&\leq& C' e^{- 2 t_0} \, .
\end{eqnarray*}

The case where $\lim_{t \rightarrow - \infty} = \infty$ and $\lim_{t \rightarrow \infty} \in \partial H$ is done similarly. \\

In the second case the geodesic connects two points in $\partial H$ . Then without loss of generality we may assume that $\Delta'\rund{R \bw_t, R \bw_t}$ is maximal for $t = 0$ . So if
$\Delta'\rund{\bw_0, \bw_0} < 2 e^{2 t_0}$ we have $g' A \cap N A_{> t_0} K = \emptyset$ . Otherwise by lemma \ref{geodesic} (ii) there exists $T \geq 0$ such
that $\Delta'\rund{\bw_T, \bw_T} = \Delta'\rund{\bw_{- T}, \bw_{- T}} = 2 e^{2 t_0}$~, and since $\Delta'\rund{\bw_T, \bw_T} \leq \frac{4}{e^{2 \abs{T}} } \Delta'\rund{\bw_0, \bw_0}$ we see that

\[
T \leq \frac{1}{2} \log\rund{2 \Delta'\rund{\bw_T, \bw_T}} - t_0 \, .
\]

So $g' a_T , g' a_{- T} \in N a_{t_0} K$ and $g' a_t \in N A_{> t_0} K$ if and only if $t \in \, ] \, - T, T \,[$ . Let $t \in \, ] \, - T, T \, [$ and assume $t \geq 0$ first. Then

\[
F\rund{g g' a_t} = F\rund{g g' a_T} - \int_{t}^{T} \widetilde f\rund{g g' a_\tau} d \tau \, ,
\]

and so

\[
\abs{F\rund{g g' a_t}} \leq \norm{F\rund{g \diamondsuit}}_{\infty, N a_{t_0} K} + \int_0^T \abs{\widetilde f\rund{g g' a_\tau}} d \tau  \, .
\]

By step II and lemma \ref{geodesic} (ii) now

\begin{eqnarray*}
\int_{0}^{T} \abs{\widetilde f\rund{g g' a_\tau}} d \tau &\leq& C' \int_{0}^{T} \frac{d \tau}{\Delta'\rund{\bw_\tau, \bw_\tau}} \\
&\leq& \frac{C'}{\Delta'\rund{\bw_0, \bw_0}} \int_{0}^{T} e^{2 \tau} d \tau \\
&\leq& \frac{C'}{2 \Delta'\rund{\bw_0, \bw_0}} e^{2 T} \\
&\leq& 2 C' e^{- 2 t_0} \, .
\end{eqnarray*}

The case $t \leq 0$ is done similarly. $\Box$ \\

(ii) Since on one hand
$\partial_\tau F \rund{\diamondsuit a_\tau} |_{\tau = 0} = \widetilde f \in L^2\rund{\Gamma \backslash G} \otimes \bigwedge\rund{\cz^r}$ and on the other
hand $\vol \rund{\Gamma \backslash G} < \infty$ it
suffices to show that $\xi F$ is bounded for all $\alpha \in \Phi \setminus \{0\}$ and $\xi \in \g^{\alpha}$~. So let $\alpha \in \Phi \setminus \{0\}$ and $\xi \in \g^{\alpha}$ . First assume
$\alpha > 0$~, which clearly implies $\alpha \geq 1$ and $\xi \in T^-$ . So there exists a continuously differentiable curve $\rho: I \rightarrow G$ contained in the $T^-$-leaf containing $1$ such that
$0 \in I$ , $\rho(0) = 1$ and $\left. \partial_t \rho(t) \right|_{t = 0} = \xi$ . Let $g \in G$ . Then by theorem \ref{differentiable} (ii) we have

\begin{eqnarray*}
\rund{\xi F}(g) &=& \left. \partial_t F\rund{g \rho(t)} \right|_{t = 0} \\
&=& - \int_0^\infty \left. \partial_t \widetilde f \rund{g \rho(t) a_\tau} \right|_{t = 0} d \tau \\
&=& - \int_0^\infty \left. \partial_t \widetilde f \rund{g a_\tau a_{- \tau} \rho(t) a_\tau} \right|_{t = 0} d \tau \\
&=& - \int_0^\infty \rund{\rund{\Ad_{a_{- \tau}} (\xi)} \widetilde f} \rund{g a_\tau} d \tau \\
&=& - \int_0^\infty e^{- \alpha \tau} \rund{\xi \widetilde f} \rund{g a_\tau} d \tau \, ,
\end{eqnarray*}

so

\[
\abs{\rund{\xi F}(g)} \leq c \norm{\xi }_2 < \infty \, ,
\]

where $c$ is the {\sc Lipschitz} constant of $\widetilde f$ . The case $\alpha < 0$ is done similarly. $\Box$ \\

Therefore by the {\sc Fourier} decomposition described above we have

\[
F = \sum_{I \in \wp(r) \, , \, \abs{I} = \rho} \phantom{1} \sum_{\nu \in \zz} F_{I \nu} \bEta^I \,  ,
\]

where $F_{I \nu} \in H_\nu$ for all $I \in \wp(r)$ , $\abs{I} = \rho$ , and $\nu \in \zz$ . $\D = \D^+ + \D^-$ , and a simple calculation shows that $\D^+$ and $\D^- \in \rz \D \oplus \g \cap \rund{T^+ \oplus T^-}$ , and
so  $\D^+ F , \D^- F \in L^2\rund{\Gamma \backslash G} \otimes \bigwedge\rund{\cz^r}$ by lemma \ref{L^2} (ii) . So we get the {\sc Fourier} decomposition of $\widetilde f$ as

\[
\widetilde f = \D F = \sum_{I \in \wp(r) \, , \, \abs{I} = \rho} \phantom{1} \sum_{\nu \in \zz} \rund{\D^+ F_{I, \nu - 2} + \D^- F_{I, \nu + 2}} \bEta^I
\]

with $\D^+ F_{I, \nu - 2} + \D^- F_{I, \nu + 2} \in H_\nu$ for all $\nu \in\zz$ . But since $f \in sS_k^{\rho}(\Gamma)$ the {\sc Fourier} decomposition of $\widetilde f$ is exactly

\[
\widetilde f = \sum_{I \in \wp(r) \, , \, \abs{I} = \rho} q_I \bEta^I
\]

with $q_I \in \C^\infty(G)^\cz \cap H_{k + \rho}$ , and so for all $I \in \wp(r)$ , $\abs{I} = \rho$ , and $\nu \in \zz$

\[
\D^+ F_{I, \nu - 2} + \D^- F_{I, \nu + 2} = \left\{ \begin{array}{cl} q_I & \text{ if } \nu = k + \rho \\
                                                                                                      0 & \text{ otherwise } \end{array} \right. \,  .
\]

\begin{lemma} $F_{I, \nu} = 0$ for $I \in \wp(r)$ , $\abs{I} = \rho$ , and $\nu \geq k + \rho$ . \end{lemma}

{\it Proof:} similar to the argument of {\sc Guillemin} and {\sc Kazhdan} in \cite{GuilKaz} . Let $I \in \wp(r)$ such that $\abs{I} = \rho$ . Then by the commutation relations of $D^+$ and $D^-$ we get for all $n \in \zz$

\begin{equation}    \label{comrel}
\norm{\D^+ F_{I, n}}_2^2 = \norm{\D^- F_{I, n}}_2^2 + \nu \norm{F_{I, n}}_2^2   \,  ,
\end{equation}

and for all $n \geq k + \rho + 1$ we have $\D^+ F_{I, n - 2} + \D^- F_{I, n + 2} = 0$ and so

\[
\norm{\D^- F_{I, n + 2}}_2 = \norm{\D^+ F_{I, n - 2}}_2  \,  .
\]

Now let $\nu \geq k + \rho$ . We will prove that

\[
\norm{\D^+ F_{I, \nu + 4 l}}_2 \geq \norm{F_{I, \nu}}_2
\]

for all $l \in \nz$ by induction on $l$ :

\begin{quote}

If $l = 0$ then the inequality is clear by (\ref{comrel}) . So let us assume that the inequality is true for some $l \in \nz$ . Then again by (\ref{comrel}) we have

\[
\norm{\D^+ F_{I, \nu + 4 l + 4}}_2^2 \geq \norm{\D^- F_{I, \nu + 4 l + 4}}_2^2 = \norm{\D^+ F_{I, \nu + 4 l}}_2^2 \geq \norm{F_{I, \nu}}_2^2  \, .
\]

\end{quote}

On the other hand $\D^+ F_I \in L^2\rund{\Gamma \backslash G}$ by lemma \ref{L^2} and so $\norm{\D^+ F_{I, n}}_2 \leadsto 0$ for $n \leadsto \infty$ . This implies
$F_\nu = 0$ . $\Box$ \\

So for all $I \in \wp(r)$ , $\abs{I} = \rho$ , we obtain $\D^+ F_{I, k + \rho - 2} = q_I$ and finally $\D^- q_I = 0$ by lemma \ref{holomorphic} , since $f \in \O(\B)$ , so

\[
\norm{q_I}_2^2 = \rund{q_I , \D^+ F_{I, k + \rho - 2}} = - \rund{\D^- q_I , F_{I, k + \rho - 2}} = 0 \, ,
\]

and so $\widetilde f = 0$ , which completes the proof of our main theorem. $\Box$

\section{Computation of the $\varphi_{\gamma_0, I, m}$ }

Fix a regular loxodromic $\gamma_0 \in \Gamma$ , $g \in G$ , $t_0 > 0$ and $w_0 \in M$ such that \\
$E_0 := E_{w_0}$ is diagonal and $\gamma_0 = g a_{t_0} w_0 g^{- 1} \in g A M g^{- 1}$ . Let \\
$D \in \rz^{r \times r}$ be diagonal such that $\exp(2 \pi i D) = E_0$ and $\chi \in \rz$ such that $j(w_0) = e^{2 \pi i \chi}$ . Now we will compute $\varphi_{\gamma_0 , I, m} \in sS_k(\Gamma)$ , $I \in \wp(r)$ , \\
$m \in \frac{1}{t_0}\rund{\zz - \rund{k + \abs{I}} \chi - \tr_I D}$ , as a relative {\sc Poincaré} series with respect to
$\Gamma_0 := \spitz{\gamma_0} \sqsubset \Gamma$ . Hereby again '$\equiv$' means equality up to a constant $\not= 0$ not necessarily independent of $\gamma_0$ , $I$ and $m$ .

\begin{theorem} Let $I \in \wp(r)$ and $k \geq 2 n + 1 - \abs{I}$ . Then for all \\
$m \in \frac{1}{t_0}\rund{\zz - \rund{k + \abs{I}} \chi - \tr_I D}$

\item[(i)]

\[
\varphi_{\gamma_0 , I, m} \equiv \sum_{\gamma \in \Gamma_0 \backslash \Gamma} q |_\gamma \in sS_k^{\rund{\abs{I}} }(\Gamma)  \, ,
\]

where

\begin{eqnarray*}
&& q := \int_{- \infty}^{\infty} e^{2 \pi i m t} \Delta\rund{\diamondsuit, g a_t \b0}^{- k - \abs{I}} \overline{j \rund{g a_t , \b0}}^{k + \abs{I}} d t \rund{E_g^{- 1} \bzeta}^I \\
&& \phantom{123456789012345678901234567890123} \in sM_k^{\rund{\abs{I}} }\rund{\Gamma_0} \cap L_k^1\rund{\Gamma_0 \backslash \B}  \, .
\end{eqnarray*}

\item[(ii)] For all $\bz \in B$ we have

\[
q \rund{\bz} \equiv \rund{\Delta \rund{\bz, \bX^+} \Delta \rund{\bz, \bX^-}}^{- \frac{k+ \abs{I}}{2}}  \rund{\frac{1 + v_1}{1 - v_1}}^{\pi i m} \rund{E_g^{- 1} \bzeta}^I \,  ,
\]

where

\[
\bX^+ := g \rund{\begin{array}{c}
1 \\
0 \\
\vdots \\
0
\end{array}}  \text{  and  } \bX^- := g \rund{\begin{array}{c}
- 1 \\
0 \\
\vdots \\
0
\end{array}}
\]

are the two fixpoints of $\gamma_0$ in $\partial B$ , and

\[
\bv := g^{- 1} \bz \in B \subset \cz^p  \, .
\]

\end{theorem}

{\it Proof:}  Let $\rho := \abs{I}$ . \\

(i) Let $f \in sS_k^{(\rho)}(\Gamma)$ , and define \\
$h = \sum_{J \in \wp(r) \, , \, \abs{J} = \rho} h_J \bEta^J \in \C^\infty \rund{\rz \times M}^\cz \otimes \bigwedge\rund{\cz^r}$ , all $h_J \in \C^\infty \rund{\rz \times M}^\cz$~, and  $b_{I, m} \in \cz$ , $m \in \frac{1}{t_0} \rund{\zz - \rund{k + \abs{I}} \chi - \tr_I D}$ , as
in theorem \ref{Fourier} . Then by standard {\sc Fourier} theory and lemma \ref{reproducing} we have

\begin{eqnarray*}
b_{I, m} &\equiv& \int_0^{t_0} e^{- 2 \pi i m t} h_I(t, 1) d t \\
&\equiv& \int_0^{t_0} e^{- 2 \pi i m t} \rund{\Delta \rund{\diamondsuit , g a_t \b0}^{- k - \rho} \rund{E_g^{- 1} \bzeta}^I , f} j \rund{g a_t , \b0}^{k + \rho} d t \\
&=& \int_0^{t_0} e^{- 2 \pi i m t} \int_G \spitz{\widetilde f, {\rund{\Delta \rund{\diamondsuit , g a_t \b0}^{- k - \rho} \rund{E_g^{- 1} \bzeta}^I}^\sim}} \times \\
&& \phantom{1} \times j \rund{g a_t , \b0}^{k + \rho} d t \,   .
\end{eqnarray*}

Since by {\sc Satake}'s theorem, theorem \ref{Satake} , $\widetilde f \in L^\infty (G) \otimes \bigwedge\rund{\cz^r}$ , and

\begin{eqnarray*}
&& \int_0^{t_0} \int_G \abs{{\rund{\Delta \rund{\diamondsuit , g a_t \b0}^{- k - \rho} \rund{E_g^{- 1} \bzeta}^I}^\sim} j \rund{g a_t , \b0}^{k + \rho} } d t \\
&& \phantom{12} = \int_0^{t_0} \int_G \abs{{\rund{\Delta \rund{\diamondsuit , \b0}^{- k - \rho} \bzeta^I}^\sim} \rund{\rund{g a_t}^{- 1} \diamondsuit}} d t \\
&& \phantom{12} \equiv \int_G \abs{\widetilde{\bzeta^I}} \\
&& \phantom{12} = \int_G \abs{j\rund{\diamondsuit, \b0}^{k + \rho} } \\
&& \phantom{12} \equiv \int_B \Delta\rund{\bZ, \bZ}^{\frac{k + \rho}{2} - (p + 1)} dV_{\Leb} < \infty   \,      ,
\end{eqnarray*}

by {\sc Tonelli}'s and {\sc Fubini}'s theorem we can interchange the order of integration:

\begin{eqnarray*}
b_{I, m} &\equiv& \int_G \spitz{\widetilde f, \int_0^{t_0} e^{2 \pi i m t} {\rund{\Delta \rund{\diamondsuit , g a_t \b0}^{- k - \rho} \rund{E_g^{- 1} \bzeta}^I}^\sim} \overline{j \rund{g a_t , \b0}}^{k + \rho} d t}     \\
&=& \rund{\int_0^{t_0} e^{2 \pi i m t} \Delta \rund{\diamondsuit , g a_t \b0}^{- k - \rho} \overline{j \rund{g a_t , \b0}}^{k + \rho} d t \rund{E_g^{- 1} \bzeta}^I, f}     \\
&=& (q, f)_{\Gamma_0} \,  ,
\end{eqnarray*}

where

\[
{\rund{\int_0^{t_0} e^{2 \pi i m t} {\Delta\rund{\diamondsuit , g a_t \b0}^{- k - \rho} \overline{j \rund{g a_t , \b0}}^{k + \rho} d t \rund{E_g^{- 1} \bzeta}^I}}^\sim} \in L^1 (G)
\otimes \bigwedge\rund{\cz^r} \, ,
\]

\[
\int_0^{t_0} e^{2 \pi i m t} \Delta \rund{\diamondsuit , g a_t \b0}^{- k - \rho} \overline{j \rund{g a_t , \b0}}^{k + \rho} d t \rund{E_g^{- 1} \bzeta}^I \in \O (\B)
\]

since $\Delta \rund{\diamondsuit , \bw} \in \O(B)$ for all $\bw \in B$ and the convergence of the integral is compact, and so by lemma \ref{Poincare}

\begin{eqnarray*}
q &:=& \sum_{\gamma' \in \Gamma_0} \left. \int_0^{t_0} e^{2 \pi i m t} \Delta \rund{\diamondsuit , g a_t \b0}^{- k - \rho} \overline{j \rund{g a_t , \b0}}^{k + \rho} d t
\rund{E_g^{- 1} \bzeta}^I \right|_{\gamma'} \\
&& \phantom{1} \in sM_k\rund{\Gamma_0} \cap L_k^1\rund{\Gamma_0 \backslash \B} \,  .
\end{eqnarray*}

Clearly

\begin{eqnarray*}
&& \left.\Delta\rund{\diamondsuit, g a_t \b0}^{- k - \rho} \rund{E_g^{- 1} \bzeta}^I \right|_{\gamma_0} \\
&& \phantom{12} = \Delta\rund{\gamma_0 \diamondsuit, g a_t \b0}^{- k - \rho} \rund{E_0 E_g^{- 1} \bzeta}^I j\rund{\gamma_0, \diamondsuit}^{k + \rho} \\
&& \phantom{12} = \Delta\rund{\diamondsuit, \gamma_0^{- 1} g a_t \b0}^{- k - \rho} \rund{E_0 E_g^{- 1} \bzeta}^I \overline{j\rund{\gamma_0^{- 1}, g a_t \b0}}^{k + \rho}   \, ,
\end{eqnarray*}

so for all $\bz \in B$ we can compute $q\rund{\bz}$ as

\begin{eqnarray*}
q\rund{\bz} &=& \sum_{\nu \in \zz} \left. \int_0^{t_0} e^{2 \pi i m t} \Delta \rund{\diamondsuit , g a_t \b0}^{- k - \rho} \rund{E_g^{- 1} \bzeta}^I
\overline{j \rund{g a_t , \b0}}^{k + \rho} d t \right|_{\gamma_0^\nu} \rund{\bz} \\
&=& \sum_{\nu \in \zz} \int_0^{t_0} e^{2 \pi i m t} \Delta\rund{\bz, \gamma_0^{- \nu} g a_t \b0}^{- k - \rho} \rund{E_0^\nu E_g^{- 1} \bzeta}^I \times \\
&& \phantom{1} \times \overline{j \rund{\gamma_0^{- \nu} g a_t , \b0}}^{k + \rho} d t  \\
&=& \sum_{\nu \in \zz} \int_0^{t_0} e^{2 \pi i m t} \Delta\rund{\bz, g a_{t - \nu t_0} \b0}^{- k - \rho} \rund{E_g^{- 1} \bzeta}^I e^{2 \pi i \nu \tr_I D} \times \\
&& \phantom{1} \times \overline{j \rund{g a_{t - \nu t_0} , \b0}}^{k + \rho} e^{2 \pi i \nu \rund{k + \rho} \chi} d t  \\
&=& \sum_{\nu \in \zz} \int_0^{t_0} e^{2 \pi i m \rund{t - \nu t_0}} \Delta\rund{\bz, g a_{t - \nu t_0} \b0}^{- k - \rho} \overline{j \rund{g a_{t - \nu t_0} , \b0}}^{k + \rho} d t \times \\
&& \phantom{1} \times \rund{E_g^{- 1} \bzeta}^I  \\
&=& \int_{- \infty}^{\infty} e^{2 \pi i m t} \Delta\rund{\bz, g a_t \b0}^{- k - \rho} \overline{j \rund{g a_t , \b0}}^{k + \rho} d t \rund{E_g^{- 1} \bzeta}^I  \, .
\end{eqnarray*}

Again by lemma \ref{Poincare} we see that $\sum_{\gamma \in \Gamma_0 \backslash \Gamma} q |_\gamma \in sM_k^{\rund{\rho} }(\Gamma) \cap L_k^1\rund{\Gamma \backslash \B}$ , and so by {\sc Satake}'s theorem, theorem \ref{Satake} , it is even
an element of $sS_k^{\rund{\rho} }(\Gamma)$ , such that

\[
b_{I, m} \equiv \rund{\sum_{\gamma \in \Gamma_0 \backslash \Gamma} q |_\gamma, f}_\Gamma  \, ,
\]

and so we conclude that $\varphi_{\gamma_0 , I, m} \equiv \sum_{\gamma \in \Gamma_0 \backslash \Gamma} q |_\gamma$ . $\Box$ \\

(ii) \begin{eqnarray*}
&& \int_{- \infty}^{\infty} e^{2 \pi i m t} \Delta\rund{\bz, g a_t \b0}^{- k - \rho} \overline{j \rund{g a_t , \b0}}^{k + \rho} d t \\
&& \phantom{12} = j\rund{g^{- 1} , \bz}^{k + \rho} \int_{- \infty}^{\infty} e^{2 \pi i m t} \Delta\rund{g^{- 1} \bz, a_t \b0}^{- k - \rho} \overline{j \rund{a_t , \b0}}^{k + \rho} d t \\
&& \phantom{12} = j\rund{g^{- 1} , \bz}^{k + \rho} \int_{- \infty}^{\infty} e^{2 \pi i m t} \rund{1 - v_1 \tanh t}^{- k - \rho} \frac{1}{\rund{\cosh t}^{k + \rho} } d t \\
&& \phantom{12} = j\rund{g^{- 1} , \bz}^{k + \rho} \int_{- \infty}^\infty \frac{e^{2 \pi i m t}}{\rund{\cosh t - v_1 \sinh t}^{k + \rho} } d t \\
&& \phantom{12} \equiv j\rund{g^{- 1} , \bz}^{k + \rho} \frac{1}{\rund{1 - v_1^2}^{\frac{k + \rho}{2}} } \rund{\frac{1 + v_1}{1 - v_1}}^{\pi i m} \\
&& \phantom{12} = j\rund{g^{- 1} , \bz}^{k + \rho} \rund{ \rund{1 - v_1} \rund{1 + v_1}}^{- \frac{k + \rho}{2}} \rund{\frac{1 + v_1}{1 - v_1}}^{\pi i m} \\
&& \phantom{12} \equiv \rund{\Delta \rund{\bz, \bX^+} \Delta \rund{\bz, \bX^-}}^{- \frac{k+ \rho}{2}}  \rund{\frac{1 + v_1}{1 - v_1}}^{\pi i m} \, . \, \Box
\end{eqnarray*}


\begin{thebibliography}{99}

\bibitem{Baily} {\sc Baily}, W. L. Jr.: Introductory lectures on Automorphic forms. Princeton University Press 1973.
\bibitem{Borth} {\sc Borthwick}, D., {\sc Klimek}, S., {\sc Lesniewski}, A. and {\sc Rinaldi}, M.: Matrix Cartan Superdomains, Super Toeplitz Operators, and Quantization. Journal
of Functional Analysis {\bf 127} (1995), 456~- 510.
\bibitem{Const} {\sc Constantinescu}, F. and {\sc de Groote}, H. F. : Geometrische und algebraische Methoden der Physik: Supermannigfaltigkeiten und Virasoro-Algebren. Teubner Verlag
Stuttgart 1994.
\bibitem{Foth} {\sc Foth}, T. and {\sc Katok}, S.: Spanning sets for automorphic forms and dynamics of the frame flow on complex hyperbolic spaces. Ergod. Th. \& Dynam. Sys. (2001),
{\bf 21}, 1071 - 1099.
\bibitem{GarlRagh} {\sc Garland}, H. and {\sc Raghunathan}, M. S.: Fundamental domains in ($\rz$-)rank 1 semisimple Lie groups. Ann. Math. {\bf 92} (2) (1970), 279 - 326.
\bibitem{GuilKaz} {\sc Guillemin}, V. and {\sc Kazhdan}, D.: Some Inverse Spectral Results for Negatively Curved 2-Manifolds, Topology {\bf 19} (1980), 301 - 312.
\bibitem{KatHas} {\sc Katok}, A. and {\sc Hasselblat}, B.: Introduction to the Modern Theory of Dynamical Systems. Cambridge University Press 1995.
\bibitem{Katok} {\sc Katok}, S.: Livshitz theorem for the unitary frame flow. Ergod. Th. \& Dynam. Sys. (2004), {\bf 24}, 127 - 140.
\bibitem{Katz} {\sc Katznelson}, Y.: An introduction to Harmonic Analysis, third edition. Cambridge University Press 2004.
\bibitem{Kne} {\sc Knevel}, R.: A {\sc Satake} type theorem for Super Automorphic forms. 2007 , to appear in the Journal of Lie Theory.
\bibitem{KneBuch} {\sc Knevel}, R.: Cusp forms, Spanning sets, and Super Symmetry, A New Geometric Approach to the Higher Rank and the Super Case. VDM Saarbrücken 2008.
\bibitem{Up} {\sc Upmeier}, H.: Toeplitz Operators and Index Theory in Several Complex Variables. Birkhäuser 1996.
\bibitem{Zim} {\sc Zimmer}, R. J.: Ergodic Theory and Semisimple Groups. Monographs in Mathematiks, 81. Birkhäuser 1984.

\end{thebibliography}
\end{document}